\documentclass[pdflatex,sn-mathphys-num]{sn-jnl}

\usepackage{graphicx}
\usepackage{multirow}
\usepackage{amsmath,amssymb,amsfonts}
\usepackage{amsthm}
\usepackage[title]{appendix}
\usepackage{xcolor}
\usepackage{textcomp}
\usepackage{manyfoot}
\usepackage{booktabs}
\usepackage{tikz}
\usepackage{pgfplots}
\usepgfplotslibrary{groupplots}
\pgfplotsset{compat=1.18}
\usetikzlibrary{arrows.meta,positioning,calc,decorations.pathreplacing}
\hypersetup{bookmarksdepth=2,
  pdftitle={The Arc Threshold of the Two-Row Ising Torus: Exact Factorization, an Algebraic Critical Coupling, Uniformity in the Circumference, and the Rate of Convergence},
  pdfauthor={Marco Mandap}}

\theoremstyle{thmstyleone}
\newtheorem{theorem}{Theorem}[section]
\newtheorem{proposition}[theorem]{Proposition}
\newtheorem{lemma}[theorem]{Lemma}
\newtheorem{corollary}[theorem]{Corollary}

\theoremstyle{thmstyletwo}

\newtheorem{remark}{Remark}[section]

\theoremstyle{thmstylethree}
\newtheorem{definition}{Definition}[section]
\newtheorem{assumption}{Assumption}[section]

\newcommand{\R}{\mathbb{R}}
\newcommand{\Z}{\mathbb{Z}}
\newcommand{\C}{\mathbb{C}}

\DeclareMathOperator{\Rea}{Re}
\DeclareMathOperator{\Ima}{Im}
\DeclareMathOperator{\Tr}{Tr}
\DeclareMathOperator{\disc}{disc}

\AtBeginDocument{\raggedbottom\hbadness=10000\vbadness=10000}

\begin{document}

\title[Arc Threshold of the Two-Row Ising Torus]{The Arc Threshold of the Two-Row Ising Torus: Exact Factorization, an Algebraic Critical Coupling, Uniformity in the Circumference, and the Rate of Convergence}

\author*[1]{\fnm{Marco} \sur{Mandap}}\email{marco.mandap@bulsu.edu.ph}

\affil*[1]{\orgdiv{Mathematics Department, College of Science}, \orgname{Bulacan State University}, \orgaddress{\city{City of Malolos}, \state{Bulacan}, \postcode{3000}, \country{Philippines}}}

\abstract{The Lee--Yang circle theorem places every zero of a ferromagnetic Ising partition polynomial on the unit circle. The arc condition, that every zero lie in the closed left half-plane, turns the circle theorem into a product formula for the characteristic function of the occupation count and makes it usable for local limit theorems. For the two-row periodic Ising torus we give the exact reflection factorization of the row transfer matrix, reduce the crossing at $\zeta=i$ to a real cubic, and identify the critical coupling $a_c^*=0.8637850566870\ldots$ as the unique root in $(0,1)$ of an explicit octic. We prove that every $a>a_c^*$ satisfies the arc condition simultaneously for all circumferences, that no smaller coupling does, and that the finite thresholds approach $a_c^*$ at the rate $C/n^2$ with an explicit algebraic constant $C$. The factorization and the $2\times2$ threshold are formalized in Lean~4 against Mathlib; the $\zeta=i$ reduction, uniformity and rate are analysis-level.}



\maketitle

\section{Introduction}\label{sec:intro}

\subsection{The circle, the arc, and the local limit theorem}\label{sec:circle}

Let $V$ be a finite site set and let
\[
  Q(\zeta)=\sum_{S\subseteq V}c(S)\zeta^{\lvert S\rvert}
\]
be the partition polynomial of a ferromagnetic Ising system, with $c(S)\ge0$. The Lee--Yang circle theorem, in the form obtained by Asano contraction and extended by Ruelle and by Lieb and Sokal \cite{YangLee1952,LeeYang1952,Asano1970,Ruelle1971,LiebSokal1981}, places every zero of $Q$ on the unit circle. The theorem is rigid, and the rigidity is useful: a zero set contained in a circle is described by one real parameter per conjugate pair.

For a local limit theorem, however, the circle is not enough. The characteristic function of the occupation count is
\[
  \varphi(t)=\frac{Q(e^{it})}{Q(1)},
\]
and if the zeros are $\zeta_1,\dots,\zeta_N$, then
\begin{equation}\label{eq:factor}
  \varphi(t)=\prod_{j=1}^{N}\frac{-\zeta_j+e^{it}}{-\zeta_j+1}.
\end{equation}
Each factor is the characteristic function of a three-point law only when $-\zeta_j\ge0$, that is, when $\Rea\zeta_j\le0$. A zero in the right half-plane produces a negative multiplier; the magnitude estimate needed on the complement of a neighbourhood of the origin then disappears. On the left arc the geometry is rigid instead: a real zero must be $-1$ and contributes a fair Bernoulli factor, while a conjugate pair contributes one two-factor block with parameters
\[
  A=-\Rea\zeta\ge0,\qquad B=-\Ima\zeta,\qquad A^2+B^2=1.
\]
This is the structure used by Lebowitz, Pittel, Ruelle and Speer \cite{LPRS2016} to obtain local limit theorems with explicit constants.

\begin{definition}[Arc condition]\label{def:arc}
A partition polynomial $Q$ satisfies the \emph{arc condition} if every zero $\zeta$ of $Q$ satisfies $\Rea\zeta\le0$.
\end{definition}

For a one-parameter family of Ising systems the arc condition becomes a question about the coupling. The zeros move continuously and remain on the circle, so the boundary of the admissible set is the set of couplings at which a zero passes through $\zeta=\pm i$. In principle this is one polynomial equation. In practice the polynomial has degree $\lvert V\rvert$, and the equation is useful only when the transfer matrix in the short direction is small.

\subsection{The two-row torus}\label{sec:why}

This paper treats the smallest genuinely two-dimensional case: a torus of height two and circumference $n+1$. There are two reasons.

First, the short-direction transfer matrix is $4\times4$, but the row interchange commutes with it. The interchange splits the space into a one-dimensional antisymmetric line and a three-dimensional symmetric block. The antisymmetric eigenvalue is the monomial $\zeta a^2(1-a^2)$; the symmetric block has an explicit cubic. The partition polynomial is therefore a sum of the $(n+1)$-st powers of four explicit eigenvalues, valid for every circumference.

Second, the family admits a threshold that is uniform in the circumference. The finite thresholds $a_c(n)$ do not drift to $1$; their supremum is an explicit algebraic number, and the approach is of order $n^{-2}$. This is the property that makes the Lee--Yang route usable for a sequence of tori at a fixed coupling.

\subsection{Results}\label{sec:results}

Let $\Lambda_n=\{0,1\}\times\mathbb Z/(n+1)\mathbb Z$ be the two-row torus of circumference $n+1$, let $a\in[0,1]$ be the coupling, and let $Q_n(\zeta,a)$ be its partition polynomial; the conventions are fixed in \S\ref{sec:model}.

\begin{itemize}
\item[\textbf{(R1)}] \emph{Trace identity.} For every height $m$ and circumference $n$, the partition function of the $(m+1)\times(n+1)$ torus is $\Tr(\mathcal T^{n+1})$, where $\mathcal T$ is the row transfer matrix on the $2^{m+1}$ column configurations; at $m=1$ it is the $4\times4$ matrix of \S\ref{sec:matrix}. (Kernel-proved.)

\item[\textbf{(R2)}] \emph{Reflection factorization.} The transfer matrix commutes with the row interchange. On the antisymmetric line it is the scalar $\lambda_4=\zeta a^2(1-a^2)$; on the symmetric subspace it is the explicit $3\times3$ matrix $\mathcal S$. Hence
\[
  Q_n(\zeta,a)=\Tr(\mathcal S^{n+1})+\bigl(\zeta a^2(1-a^2)\bigr)^{n+1}
\]
for every $n\ge0$. (Kernel-proved.)

\item[\textbf{(R3)}] \emph{The cubic and its recurrence.} $\mathcal S$ satisfies the cubic
\[
  X^3-\tau X^2+\sigma X-\delta,
\]
with $\tau,\sigma,\delta$ given in \S\ref{sec:factor}; its traces obey the cubic recurrence. The partition polynomial is the sum of the $(n+1)$-st powers of the three roots and the antisymmetric monomial. (Kernel-proved.)

\item[\textbf{(R4)}] \emph{The $2\times2$ threshold.} For the $2\times2$ torus, $a^4\ge\sqrt2-1$ implies the arc condition, and the constant is exact: at $a^4=\sqrt2-1$ the polynomial vanishes at $\zeta=i$. (Kernel-proved; sharpness is an elementary substitution.)

\item[\textbf{(R5)}] \emph{The critical coupling.} At $\zeta=i$ the cubic becomes $X^3-i\alpha X^2-\beta X+i\gamma$ with real $\alpha,\beta,\gamma$, so its roots are $i$ times the roots of the real cubic $R_a(y)=y^3-\alpha y^2+\beta y-\gamma$. Thus
\[
  Q_n(i,a)=i^{n+1}F_n(a),\qquad
  F_n(a)=p_{n+1}(a)+d(a)^{n+1},
\]
where $p_k$ is the $k$-th power sum of the roots of $R_a$ and $d=a^2(1-a^2)$. The discriminant of $R_a$ factors as $4(a^2-1)^2P(a^2)$ with $P$ an explicit octic, and the unique root of $P$ in $(0,1)$ is
\[
  a_c^*=0.863785056687024029921664\ldots .
\]
(Analysis-level.)

\item[\textbf{(R6)}] \emph{Uniformity and optimality.} For $a>a_c^*$ the arc condition holds for $Q_n$ simultaneously for every $n$; for $a<a_c^*$ it fails for some $n$. Hence $\sup_n a_c(n)=a_c^*$. (Analysis-level; the zero-count input is Rouch\'e's theorem in Appendix~\ref{app:rouche}, and the dominance lemma is proved in \S\ref{sec:uniform}.)

\item[\textbf{(R7)}] \emph{The rate.} The finite thresholds satisfy
\[
  a_c(n)=a_c^*-\frac{C}{n^2}+\frac{2C}{n^3}-\frac{3C+D}{n^4}+O(n^{-5}),
\]
with explicit algebraic
\[
  C=\frac{\pi^2\delta_0^4u^2}{\disc'(a_c^*)}=0.1914850\ldots,
  \qquad D=0.130886\ldots .
\]
(Analysis-level; Appendix~\ref{app:rate}.)
\end{itemize}

The finite thresholds are recorded in \S\ref{sec:finite}: $a_c(1)=(\sqrt2-1)^{1/4}=0.802243\ldots$, $a_c(10)=0.862193\ldots$, $a_c(40)=0.863671\ldots$, increasing toward $a_c^*$ with a gap of order $n^{-2}$.

\subsection{Method}\label{sec:method}

Three algebraic facts organize the paper.

The first is the trace identity. The partition function of a height-$h$ torus is the trace of the $(n+1)$-st power of a $2^h\times2^h$ row transfer matrix. At height two this reduces a sum over $2^{2(n+1)}$ spin configurations to a matrix power.

The second is the row interchange. A column's weight depends only on the multiset of its two spins and on the unordered pair of configurations, so the interchange commutes with the transfer matrix. The one-dimensional antisymmetric line has a monomial eigenvalue; the three-dimensional symmetric block has a cubic.

The third is the evaluation at $\zeta=i$. On the unit circle this is the only point where the right half-plane meets the circle, so it is the only place where a zero can cross as $a$ varies. The factorization reduces the crossing to the real cubic $R_a$; positivity of its roots becomes a discriminant condition, and the discriminant produces the octic $P$.

Uniformity and the rate are then analytic consequences of the same cubic. Above $a_c^*$ the cubic has three positive roots, so every term of $F_n$ is positive; below $a_c^*$ the complex pair dominates the real root and the antisymmetric monomial. The rate comes from the implicit function theorem applied to the simple zero of $F_n$ near $a_c^*$, with the constants read off from the cubic and the octic.

\subsection{Proof standards}\label{sec:status}

The paper uses two standards and does not blur them.

\emph{Kernel-proved.} Results (R1)--(R4) are theorems in Lean~4 against Mathlib, in
\texttt{SteinGate.LeeYang.LeeYang2DTorus} and
\texttt{SteinGate.LeeYang.LeeYang2DTorusArcCondition}. Their axiom footprint is Mathlib's
(\texttt{propext}, \texttt{Classical.choice}, \texttt{Quot.sound}), and neither module contains a
\texttt{sorry}. Appendix~\ref{app:lean} lists the declarations and the commands.

\emph{Analysis-level.} Results (R5)--(R7) are proved on paper. Their inputs are the circle theorem, continuity of zeros in the coefficients, Rouch\'e's theorem, the Beraha--Kahane--Weiss equimodular limit, and the implicit function theorem. The arguments are complete at that level; no proof assistant has checked them. The numerical tables are certified computations, not proofs.

\emph{Open.} Extensions beyond height two, uniform-in-height structure, and fixed boundary conditions are outside this paper and are stated as open in \S\ref{sec:discussion}.

\subsection{Organization}\label{sec:org}

\S\ref{sec:model} fixes the model and the trace identity. \S\ref{sec:matrix} introduces the $4\times4$ matrix and the reflection decomposition. \S\ref{sec:factor} gives the factorization and the cubic. \S\ref{sec:two} treats the $2\times2$ torus. \S\ref{sec:i} evaluates at $\zeta=i$; \S\ref{sec:crit} identifies $a_c^*$; \S\ref{sec:finite} records the finite thresholds. \S\ref{sec:uniform} proves uniformity, optimality, and the exact threshold. \S\ref{sec:discussion} places the result in the literature and lists the open directions. Appendix~\ref{app:lean} documents the formal artifact; Appendix~\ref{app:rate} proves the rate.

\section{The model and the transfer matrix}\label{sec:model}

\subsection{The successor-edge torus}\label{sec:edges}

Fix $n\ge 1$ and work on the site set
\[
  \Lambda_{n}\;=\;\{0,1\}\times\Z/(n+1)\Z ,
\]
with the second coordinate taken modulo $n+1$. For a site $p=(i,j)$ write $p+e_{1}=(i,j+1)$ and $p+e_{2}=(i+1,j)$, both coordinates modulo $2$ and $n+1$ respectively.

\begin{definition}[Successor-edge torus]\label{def:torus}
The edge list $\mathcal{E}_{n}$ is the list of ordered pairs
\[
  \mathcal{E}_{n}\;=\;\bigl\{(p,p+e_{1}) : p\in\Lambda_{n}\bigr\}\;\cup\;\bigl\{(p,p+e_{2}) : p\in\Lambda_{n}\bigr\},
\]
each carrying the same coupling $a\in[0,1]$. An ordered pair $(p,q)$ \emph{crosses} a spin set $S\subseteq\Lambda_{n}$ if exactly one of $p,q$ lies in $S$. The \emph{partition function} and \emph{partition polynomial} are
\[
  Z_{n}(\zeta,a)=\sum_{S\subseteq\Lambda_{n}}a^{c(S)}\zeta^{\lvert S\rvert},
  \qquad
  Q_{n}(\zeta,a)=\sum_{S\subseteq\Lambda_{n}}a^{c(S)}\zeta^{\lvert S\rvert}\in\Z[a][\zeta],
\]
where $c(S)$ is the number of edges of $\mathcal{E}_{n}$ crossing $S$. Thus $Q_{n}(\cdot,a)$ is a monic polynomial of degree $2(n+1)$.
\end{definition}

\begin{remark}[Why the edge list is a list]\label{rem:doubling}
Because $\mathcal{E}_{n}$ is a list of ordered pairs, a bond whose two orientations both occur contributes the factor $a^{2}$ rather than $a$. This happens exactly when the successor map in the relevant direction has period two. In the vertical direction the period is always two, so every vertical bond of a two-row torus carries weight $a^{2}$; in the horizontal direction the period is $n+1$, so the horizontal bonds carry weight $a$ unless $n=1$, in which case they too carry weight $a^{2}$. Equivalently: the model has horizontal coupling $a$ and vertical coupling $a^{2}$, and for $n=1$ both couplings are $a^{2}$. This is the convention of the formal artifact and is exactly what makes the transfer matrix symmetric under the interchange of the two rows.
\end{remark}

\begin{figure}[tbp]
\centering
\begin{tikzpicture}[x=1.15cm,y=1.15cm,font=\small,>=Latex]
 \fill[blue!7,rounded corners=3pt] (1.70,-0.26) rectangle (2.30,1.41);
 \foreach \i in {0,1}{
   \foreach \j in {0,...,4}{
     \draw[black!65,line width=.55pt] (\j,\i*1.15) -- (\j+1,\i*1.15);
   }
 }
 \foreach \j in {0,...,5}{
   \draw[blue!60!black,line width=.45pt] (\j-0.032,0.055) -- (\j-0.032,1.095);
   \draw[blue!60!black,line width=.45pt] (\j+0.032,0.055) -- (\j+0.032,1.095);
 }
 \foreach \j in {0,...,5}{
   \foreach \i in {0,1}{
     \fill[black!85] (\j,\i*1.15) circle (1.9pt);
   }
 }
 \node[above,font=\small] at (3.50,0.10) {$a$};
 \node[anchor=west,font=\small,blue!60!black] at (4.11,0.575) {$a^{2}$};
 \node[anchor=south,font=\small] at (2.00,1.44) {$b\in\{0,1\}^{2}$};
 \node[anchor=west,font=\footnotesize,black!70] at (5.40,0.00) {row $0$};
 \node[anchor=west,font=\footnotesize,black!70] at (5.40,1.15) {row $1$};
 \draw[black!40,densely dotted,line width=.4pt] (-0.45,-0.30) -- (-0.45,1.44);
 \draw[black!40,densely dotted,line width=.4pt] (5.25,-0.30) -- (5.25,1.44);
 \draw[-{Latex[length=1.6mm]},black!45,dashed,line width=.4pt]
   (5.25,1.44) .. controls (5.25,2.10) and (-0.45,2.10) .. (-0.45,1.48);
 \node[font=\footnotesize,black!70] at (2.40,2.10)
   {periodic in the horizontal direction};
 \draw[->,black!60,line width=.5pt] (2.32,-0.58) -- (2.98,-0.58);
 \node[below,font=\footnotesize] at (2.65,-0.62) {$\mathcal{T}(\zeta,a)$};
 \node[font=\small] at (-0.95,1.95) {(a)};
 \foreach \x/\sa/\sb/\lab in {
     1.05/0/0/{$c_{00}$},
     2.55/1/0/{$c_{10}$},
     4.05/0/1/{$c_{01}$},
     5.55/1/1/{$c_{11}$}}{
   \draw[black!45,line width=.4pt] (\x-0.19,-2.47) rectangle (\x+0.19,-2.09);
   \draw[black!45,line width=.4pt] (\x-0.19,-2.85) rectangle (\x+0.19,-2.47);
   \ifnum\sa=1
     \fill[black!80] (\x,-2.28) circle (1.9pt);
   \else
     \draw[black!50,line width=.5pt] (\x,-2.28) circle (1.9pt);
   \fi
   \ifnum\sb=1
     \fill[black!80] (\x,-2.66) circle (1.9pt);
   \else
     \draw[black!50,line width=.5pt] (\x,-2.66) circle (1.9pt);
   \fi
   \node[below,font=\small] at (\x,-2.90) {\lab};
 }
 \node[anchor=east,font=\footnotesize,black!70] at (0.78,-2.28) {$0$};
 \node[anchor=east,font=\footnotesize,black!70] at (0.78,-2.66) {$1$};
 \draw[decorate,decoration={brace,amplitude=3pt},black!55]
   (2.28,-2.02) -- (4.32,-2.02)
   node[midway,above=4pt,font=\small] {$\rho$};
 \node[font=\small] at (-0.95,-2.47) {(b)};
\end{tikzpicture}
\caption{The lattice and the transfer step, drawn schematically for
$n+1=6$. (a)~The $2\times(n+1)$ torus. Horizontal bonds carry the weight
$a$; each vertical bond carries $a^{2}$, because both orientations of the
successor edge occur (Remark~\ref{rem:doubling}). The two dotted lines are
identified, so the horizontal direction is a cycle of length $n+1$. The
shaded band is one column configuration $b$; the arrow is the row-transfer
step to the next column. (b)~The four column configurations, one square per
row, with a filled site for spin $1$ and an open site for spin $0$; the
involution $\rho$ exchanges $c_{10}$ and $c_{01}$ and fixes $c_{00}$ and
$c_{11}$.}
\label{fig:torus}
\end{figure}
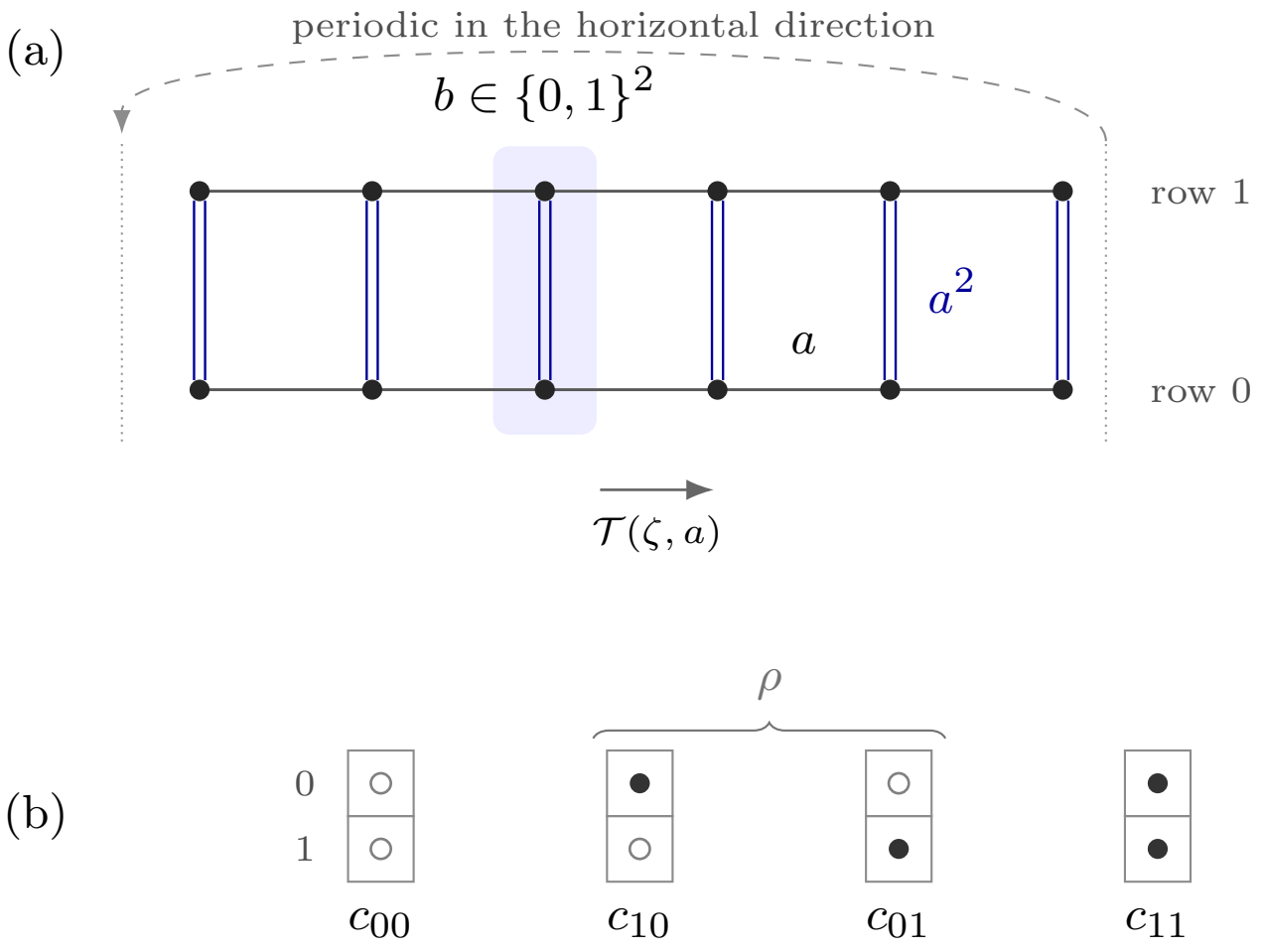

Figure~\ref{fig:torus} records both conventions: the two bond weights on the torus, and the four column configurations on which the row transfer matrix acts.

\subsection{The row transfer matrix}\label{sec:tm}

A \emph{column configuration} is a function $b:\{0,1\}\to\{0,1\}$, assigning a spin to the two sites of a column. There are four of them; we label them
\[
  c_{00}=(0,0),\qquad c_{10}=(1,0),\qquad c_{01}=(0,1),\qquad c_{11}=(1,1),
\]
the subscript recording the spins in rows $0$ and $1$ in that order.

\begin{definition}[Row transfer matrix]\label{def:tm}
For $\zeta\in\C$ and $a\in[0,1]$ let $\mathcal{T}(\zeta,a)$ be the $4\times 4$ complex matrix indexed by column configurations,
\[
  \begin{aligned}
  \mathcal{T}(\zeta,a)_{b,c}
  \;=\;
  &\Biggl(\prod_{i=0}^{1}\bigl(\zeta\,\mathbf{1}\{b_{i}=1\}\bigr)
    \bigl(\mathbf{1}\{b_{i}=b_{i+1}\}+a\,\mathbf{1}\{b_{i}\ne b_{i+1}\}\bigr)\Biggr)\\
  &\times\Biggl(\prod_{i=0}^{1}\bigl(\mathbf{1}\{b_{i}=c_{i}\}+a\,\mathbf{1}\{b_{i}\ne c_{i}\}\bigr)\Biggr),
  \end{aligned}
\]
with $b_{2}:=b_{0}$.
\end{definition}

The first product collects the within-column weight (the vertical bonds, each counted twice by Remark~\ref{rem:doubling}, and the fugacity of the occupied sites of the column); the second collects the between-column weight (the horizontal bonds joining column $b$ to column $c$).

\begin{theorem}[Trace identity]\label{thm:trace}
For every $m,n\ge 0$, every $a\in[0,1]$ and every $\zeta\in\C$, the partition function of the $(m+1)\times(n+1)$ torus equals $\Tr\bigl(\mathcal{T}_{m}(\zeta,a)^{n+1}\bigr)$, where $\mathcal{T}_{m}$ acts on the $2^{m+1}$ column configurations $\{0,\dots,m\}\to\{0,1\}$ by the natural analog of Definition~\ref{def:tm}. In particular, for $m=1$,
\[
  Z_{n}(\zeta,a)\;=\;\Tr\bigl(\mathcal{T}(\zeta,a)^{n+1}\bigr).
\]
\end{theorem}

\begin{proof}
Kernel-proved: \texttt{isingPartitionFun\_torusEdges\_eq\_trace} in \texttt{SteinGate.LeeYang.LeeYang2DTorus} for the general case, and its specialization at $m=1$. The argument expands the sum over spin sets column by column, recognizes the product over columns as a matrix product, and closes the ring using the cyclic indexing. Appendix~\ref{app:lean} records the exact statements.
\end{proof}

\begin{remark}\label{rem:traceproof}
The proof is a closed-walk count on a transfer graph whose vertices are column configurations and whose weighted adjacency matrix is $\mathcal{T}$. The trace appears because the horizontal direction is a cycle of length $n+1$; the vertical direction is a cycle of length $2$, which is the origin of the square in Remark~\ref{rem:doubling}.
\end{remark}

\subsection{The arc condition and its threshold}\label{sec:threshold}

\begin{definition}[Arc condition; threshold]\label{def:acn}
Fix $n\ge1$ and $a\in[0,1]$. Say that $Q_{n}(\cdot,a)$ satisfies the \emph{arc condition} if every zero of $Q_{n}(\cdot,a)$ satisfies $\Rea\zeta\le 0$. Write
\[
  A_{n}\;=\;\bigl\{a\in[0,1] : Q_{n}(\cdot,a)\text{ satisfies the arc condition}\bigr\},
  \qquad
  a_{c}(n)\;=\;\inf A_{n}.
\]
\end{definition}

Two remarks. First, $A_{n}$ is nonempty: at $a=1$ all weights are $1$, so $Q_{n}(\zeta,1)=(1+\zeta)^{2(n+1)}$, and continuity of the zeros gives a neighborhood of $1$. Second, the zeros lie on the circle and $\{\Rea\zeta=0\}=\{\pm i\}$ there, so a zero can enter or leave the open right half-plane only through $\pm i$; hence $A_{n}$ is a union of intervals whose endpoints are roots of $Q_{n}(\pm i,\cdot)$. In fact $A_{n}=[a_{c}(n),1]$ (\S\ref{sec:uniform}, analysis-level).

\begin{assumption}[Circle theorem]\label{ass:circle}
For every $n\ge1$, every $a\in[0,1]$ and every zero $\zeta$ of $Q_{n}(\cdot,a)$, one has $\lvert\zeta\rvert=1$.
\end{assumption}

Assumption~\ref{ass:circle} is a theorem \cite{YangLee1952,LeeYang1952,Asano1970,Ruelle1971,LiebSokal1981}, not an extra hypothesis; we isolate it because it is the one structural input of \S\ref{sec:uniform} that the formal artifact does not supply for this particular graph. The project's Lean library proves the circle theorem for the models treated in \texttt{SteinGate.LeeYang.LeeYangCircle} and, at the torus, in \texttt{torusEdges\_roots\_norm}; the general finite-graph version is the classical Asano contraction.

\section{The \texorpdfstring{$4\times 4$}{4x4} matrix and its reflection decomposition}\label{sec:matrix}

\subsection{The matrix in closed form}\label{sec:matrixclosed}

\begin{lemma}[Entries]\label{lem:entries}
With the row and column order $c_{00},c_{10},c_{01},c_{11}$,
\begin{equation}\label{eq:T}
  \mathcal{T}(\zeta,a)\;=\;
  \begin{pmatrix}
    1 & a & a & a^{2}\\[2pt]
    \zeta a^{3} & \zeta a^{2} & \zeta a^{4} & \zeta a^{3}\\[2pt]
    \zeta a^{3} & \zeta a^{4} & \zeta a^{2} & \zeta a^{3}\\[2pt]
    \zeta^{2}a^{2} & \zeta^{2}a & \zeta^{2}a & \zeta^{2}
  \end{pmatrix}.
\end{equation}
\end{lemma}

\begin{proof}
Kernel-proved entry by entry: the lemmas \texttt{T\_c00\_c00}, \dots, \texttt{T\_c11\_c11} in \texttt{SteinGate.LeeYang.LeeYang2DTorusArcCondition}. Each is a direct unfolding of Definition~\ref{def:tm} at the relevant pair of configurations.
\end{proof}

The factor $2$ that appears in the middle row and column of the symmetric block below is the sum of the two entries $\zeta a^{4}$ of \eqref{eq:T} (the pair $c_{10}\leftrightarrow c_{01}$), and the corresponding $2$ in the bottom-left corner is the sum of the two entries $\zeta^{2}a$.

\subsection{The row interchange}\label{sec:reflection}

Let $\rho$ be the involution of column configurations that interchanges the two rows,
\[
  (\rho b)_{0}=b_{1},\qquad (\rho b)_{1}=b_{0},
\]
so that $\rho$ swaps $c_{10}$ and $c_{01}$ and fixes $c_{00}$ and $c_{11}$. It acts on functions $x$ on column configurations by $(Px)(b)=x(\rho b)$, and $P^{2}=\mathrm{id}$.

\begin{lemma}[Reflection symmetry]\label{lem:reflection}
$P\,\mathcal{T}(\zeta,a)\,P=\mathcal{T}(\zeta,a)$.
\end{lemma}

\begin{proof}
Immediate from \eqref{eq:T}: the entry $(b,c)$ depends on $b$ and $c$ only through the unordered pairs $\{b_{0},b_{1}\}$ and the equality pattern of $b$ with $c$, both invariant under simultaneous interchange of rows in $b$ and in $c$. Equivalently, the weight of a column in Definition~\ref{def:tm} is unchanged by $\rho$.
\end{proof}

\begin{figure}[tbp]
\centering
\begin{tikzpicture}[font=\small,>=Stealth,
  box/.style={draw=black!60,rounded corners=2pt,fill=black!2,align=center,
    inner sep=5pt},
  ar/.style={-{Stealth[length=2mm]},line width=.6pt,draw=black!70}]
 \node[box,text width=9.4cm,minimum height=0.95cm] (T) at (0,0)
   {$\mathcal{T}(\zeta,a)$ on $\C^{4}$, basis $c_{00},c_{10},c_{01},c_{11}$};
 \node[anchor=west,font=\footnotesize,black!70] at (4.82,0) {$4=3+1$};
 \node[box,text width=6.4cm,minimum height=2.7cm] (plus) at (-2.45,-2.30)
   {$V_{+}=\{x : x(c_{10})=x(c_{01})\}$, $\dim V_{+}=3$\\[3pt]
    coordinates $(x_{0},x_{1},x_{2})$\\[3pt]
    $\mathcal{S}(\zeta,a)$: the $3\times3$ block \eqref{eq:S}\\[3pt]
    cubic $X^{3}-\tau X^{2}+\sigma X-\delta$};
 \node[box,text width=3.9cm,minimum height=2.7cm] (minus) at (3.95,-2.30)
   {$V_{-}=\C\cdot v_{-}$, $\dim V_{-}=1$\\[3pt]
    $v_{-}=(0,1,-1,0)$\\[3pt]
    $\lambda_{4}=\zeta a^{2}(1-a^{2})$};
 \draw[ar] (-2.45,-0.48) -- (-2.45,-0.95)
   node[midway,left=3pt,font=\footnotesize] {$P=+1$};
 \draw[ar] (3.95,-0.48) -- (3.95,-0.95)
   node[midway,right=3pt,font=\footnotesize] {$P=-1$};
 \draw[ar] (-2.45,-3.65) -- (-2.45,-3.95) -- (3.95,-3.95) -- (3.95,-3.65);
 \draw[ar] (0.75,-3.95) -- (0.75,-4.22);
 \node[box,text width=8.6cm,minimum height=0.80cm] (Q) at (0.75,-4.65)
   {$Q_{n}(\zeta,a)=\Tr\bigl(\mathcal{S}(\zeta,a)^{n+1}\bigr)
     +\bigl(\zeta a^{2}(1-a^{2})\bigr)^{n+1}$};
\end{tikzpicture}
\caption{The block decomposition, schematically. The row interchange splits
$\C^{4}$ into the $P$-eigenspaces $V_{+}$ and $V_{-}$; $\mathcal{T}$ respects
the splitting, acting on $V_{+}$ as the $3\times3$ block $\mathcal{S}$ of
\eqref{eq:S} and on $V_{-}$ as the scalar $\lambda_{4}$. Taking traces in the
trace identity of Theorem~\ref{thm:trace} gives the factorization
\eqref{eq:factorization}: three powers of the roots of the cubic plus one
explicit power.}
\label{fig:blocks}
\end{figure}
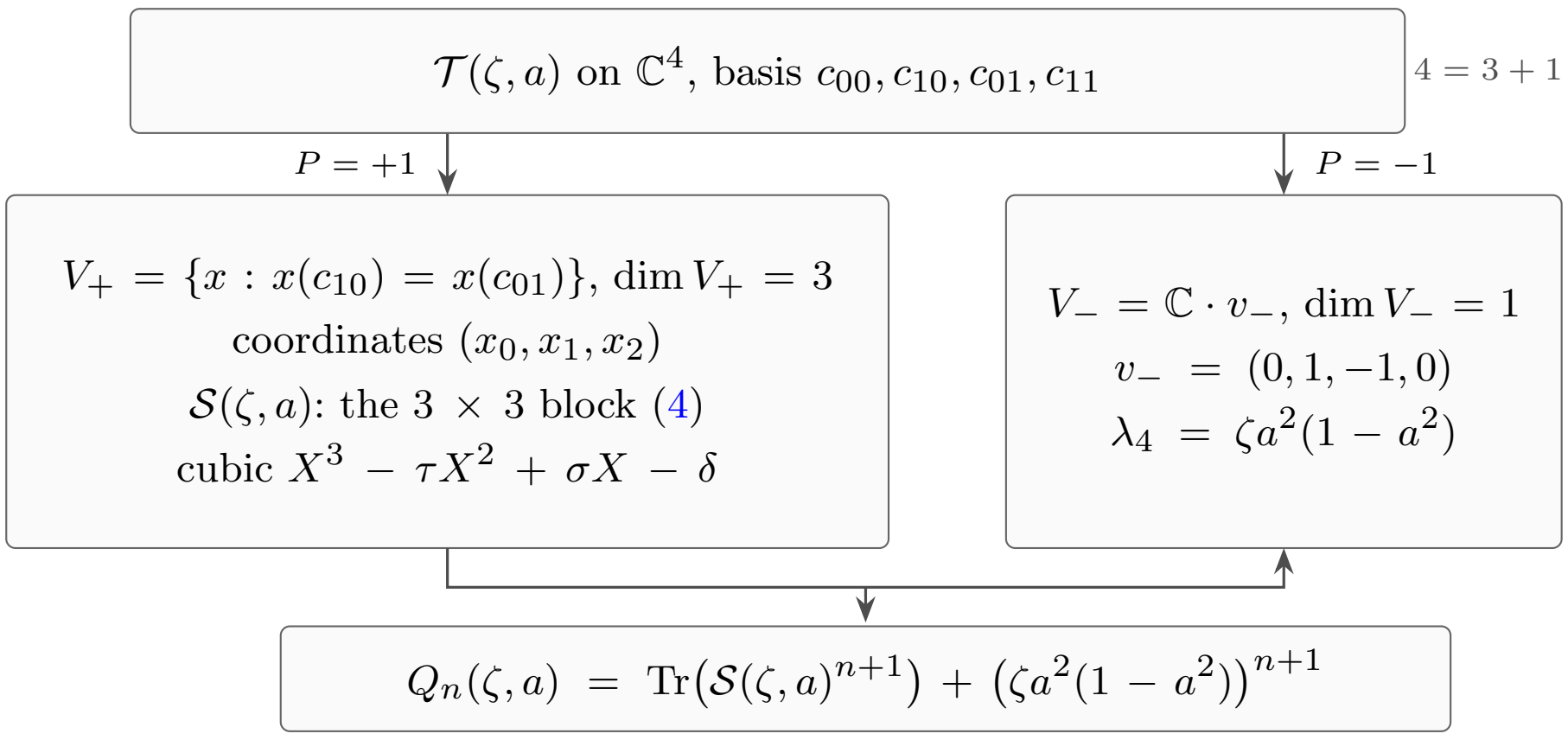

Hence $\C^{4}$ decomposes as the direct sum of the $P$-eigenspaces (Figure~\ref{fig:blocks}). The antisymmetric space is
\[
  V_{-}\;=\;\bigl\{x : x(c_{00})=x(c_{11})=0,\; x(c_{10})=-x(c_{01})\bigr\}\;=\;\C\cdot v_{-},
  \qquad v_{-}=(0,1,-1,0),
\]
one-dimensional, and the symmetric space
\[
  V_{+}\;=\;\bigl\{x : x(c_{10})=x(c_{01})\bigr\}
\]
is three-dimensional.

\begin{lemma}[Antisymmetric eigenvalue]\label{lem:antisym}
$\mathcal{T}(\zeta,a)\,v_{-}=\lambda_{4}(\zeta,a)\,v_{-}$ with
\begin{equation}\label{eq:lambda4}
  \lambda_{4}(\zeta,a)\;=\;\zeta a^{2}(1-a^{2}).
\end{equation}
\end{lemma}

\begin{proof}
Kernel-proved: \texttt{torusT\_antisym\_eigenvalue}. For the reader's convenience: the $c_{10}$-row of \eqref{eq:T} applied to $v_{-}$ gives $\zeta a^{2}\cdot 1+\zeta a^{4}\cdot(-1)=\zeta a^{2}(1-a^{2})$, and the $c_{01}$-row gives $\zeta a^{4}-\zeta a^{2}$, which is the negative of the same number; the other two rows give $a-a=0$. One checks in the same way that $v_{-}$ is the whole antisymmetric space, so no further antisymmetric eigenvalue exists.
\end{proof}

\subsection{The symmetric block}\label{sec:symblock}

On $V_{+}$ use the coordinates $x_{0}=x(c_{00})$, $x_{1}=x(c_{10})=x(c_{01})$, $x_{2}=x(c_{11})$.

\begin{definition}[Symmetric block]\label{def:block}
\begin{equation}\label{eq:S}
  \mathcal{S}(\zeta,a)\;=\;
  \begin{pmatrix}
    1 & 2a & a^{2}\\[2pt]
    \zeta a^{3} & \zeta a^{2}(1+a^{2}) & \zeta a^{3}\\[2pt]
    \zeta^{2}a^{2} & 2\zeta^{2}a & \zeta^{2}
  \end{pmatrix}.
\end{equation}
\end{definition}

\begin{lemma}[Invariance]\label{lem:invariance}
For every $x\in\C^{3}$, the vector $x\mapsto(x_{0},x_{1},x_{1},x_{2})$ is mapped by $\mathcal{T}$ to the vector with the same symmetry, and its three coordinates are $\mathcal{S}(\zeta,a)x$.
\end{lemma}

\begin{proof}
Kernel-proved: \texttt{torusT\_mulVec\_symVec}, via the expansion lemma \texttt{torusT\_sym\_mulVec\_expand}. The verification is a four-case computation on column configurations; the factor $2$ in the $(0,1)$ and $(2,1)$ entries arises because the two preimages $c_{10}$ and $c_{01}$ both contribute.
\end{proof}

\begin{theorem}[Trace splitting]\label{thm:split}
For every $k\ge0$,
\[
  \Tr\bigl(\mathcal{T}(\zeta,a)^{k}\bigr)\;=\;\Tr\bigl(\mathcal{S}(\zeta,a)^{k}\bigr)+\lambda_{4}(\zeta,a)^{k}.
\]
\end{theorem}

\begin{proof}
Kernel-proved: \texttt{trace\_pow\_torusT}. The basis $(c_{00},\,c_{10}+c_{01},\,c_{11},\,v_{-})$ diagonalizes the action of $P$, hence block-diagonalizes $\mathcal{T}$ into $\mathcal{S}$ and the scalar $\lambda_{4}$; the trace of a power is the sum of the powers of the eigenvalues, block by block.
\end{proof}

\section{The factorization, the cubic, and the recurrence}\label{sec:factor}

Combining Theorem~\ref{thm:trace} with Theorem~\ref{thm:split} gives the factorization that organizes everything else.

\begin{theorem}[Factorization]\label{thm:factor}
For every $n\ge0$, every $a\in[0,1]$ and every $\zeta\in\C$,
\begin{equation}\label{eq:factorization}
  Q_{n}(\zeta,a)\;=\;\Tr\bigl(\mathcal{S}(\zeta,a)^{n+1}\bigr)+\bigl(\zeta a^{2}(1-a^{2})\bigr)^{n+1}.
\end{equation}
\end{theorem}

\begin{proof}
Kernel-proved: \texttt{isingPartitionFun\_torusEdges\_2xn\_eq}. Apply Theorem~\ref{thm:trace} with $m=1$ and then Theorem~\ref{thm:split} with $k=n+1$.
\end{proof}

It remains to study the three eigenvalues of $\mathcal{S}$ and the explicit scalar \eqref{eq:lambda4}.

\begin{theorem}[The cubic]\label{thm:cubic}
$\mathcal{S}(\zeta,a)$ has characteristic polynomial
\begin{equation}\label{eq:cubic}
  \chi_{\mathcal{S}}(X)\;=\;X^{3}-\tau X^{2}+\sigma X-\delta,
\end{equation}
where
\begin{align}
  \tau&=\;1+\zeta a^{2}(1+a^{2})+\zeta^{2}, \label{eq:tau}\\
  \sigma&=\;\zeta a^{2}(1-a^{2})+\zeta^{2}(1-a^{4})+\zeta^{3}a^{2}(1-a^{2}), \label{eq:sigma}\\
  \delta&=\;\zeta^{3}a^{2}(1-a^{2})^{3}. \label{eq:delta}
\end{align}
\end{theorem}

\begin{proof}
Kernel-proved: \texttt{charpoly\_torusT\_sym}, together with the definition \texttt{torusT\_sym\_charPoly} and the identification lemma \texttt{torusT\_sym\_charPoly\_eq} which names the three coefficients. The proof evaluates the determinant of $X\mathbf{1}-\mathcal{S}$ for general $X$ and compares coefficients.
\end{proof}

\begin{remark}[Reading the coefficients]\label{rem:coeffs}
The coefficients have transparent meanings. $\tau=\Tr\mathcal{S}$ is the sum of the diagonal entries of \eqref{eq:S}. $\sigma$ is the sum of the three principal $2\times 2$ minors; the three terms correspond to the pairs $(0,1)$, $(0,2)$ and $(1,2)$, and the reader may check that the $(0,1)$ minor is $\zeta a^{2}(1+a^{2})-2\zeta a^{4}=\zeta a^{2}(1-a^{2})$, and similarly for the others. Finally $\delta=\det\mathcal{S}$, and it is a pleasant check that $\det\mathcal{S}\cdot\lambda_{4}=\det\mathcal{T}=\zeta^{4}a^{4}(1-a^{2})^{4}$.
\end{remark}

\begin{theorem}[Cayley--Hamilton]\label{thm:CH}
$\mathcal{S}^{3}=\tau\,\mathcal{S}^{2}-\sigma\,\mathcal{S}+\delta\,\mathbf{1}$.
\end{theorem}

\begin{proof}
Kernel-proved: \texttt{torusT\_sym\_cayley\_hamilton}, obtained from Mathlib's \texttt{Matrix.aeval\_self\_charpoly} together with Theorem~\ref{thm:cubic}.
\end{proof}

\begin{theorem}[Recurrence]\label{thm:rec}
Let $g_{k}=\Tr(\mathcal{S}^{k})$, with $g_{0}=3$. Then for every $k\ge0$,
\begin{equation}\label{eq:recurrence}
  g_{k+3}\;=\;\tau\,g_{k+2}-\sigma\,g_{k+1}+\delta\,g_{k}.
\end{equation}
\end{theorem}

\begin{proof}
Kernel-proved: \texttt{torusT\_symTrace\_recurrence}. Multiply Theorem~\ref{thm:CH} by $\mathcal{S}^{k}$ on the left, take traces, and use cyclicity.
\end{proof}

\begin{corollary}[Power-sum form]\label{cor:powersum}
Let $y_{1},y_{2},y_{3}$ be the roots of \eqref{eq:cubic}, counted with multiplicity. Then $g_{k}=y_{1}^{k}+y_{2}^{k}+y_{3}^{k}$ for every $k\ge0$, and
\begin{equation}\label{eq:Qpowersum}
  Q_{n}(\zeta,a)\;=\;y_{1}^{n+1}+y_{2}^{n+1}+y_{3}^{n+1}+\bigl(\zeta a^{2}(1-a^{2})\bigr)^{n+1}.
\end{equation}
\end{corollary}

\begin{proof}
The characteristic polynomial of $\mathcal{S}$ is \eqref{eq:cubic}, so the eigenvalues of $\mathcal{S}$ are $y_{1},y_{2},y_{3}$; the trace of $\mathcal{S}^{k}$ is the sum of the $k$-th powers of the eigenvalues. Substitute into \eqref{eq:factorization}.
\end{proof}

\begin{remark}[What is and is not claimed]\label{rem:whatnot}
Corollary~\ref{cor:powersum} is a closed form for $Q_{n}$ in terms of the roots of a single cubic. It is not a diagonalization of the full $2^{n+1}\times 2^{n+1}$ row-to-row transfer matrix of the two-dimensional model; that matrix is what one would need for the standard Onsager-type analysis in which the circumference is the small direction. Here the small direction is the height, which is fixed at two, and the closed form is uniform in the circumference. The two analyses answer different questions and neither contains the other.
\end{remark}

\section{The \texorpdfstring{$2\times2$}{2x2} torus: an exact threshold}\label{sec:two}

For $n=1$ the torus has four sites and eight oriented edges, and by Remark~\ref{rem:doubling} every one of the four bonds carries weight $a^{2}$. The partition polynomial is a reciprocal quartic.

\begin{lemma}[Explicit quartic]\label{lem:quartic}
\begin{equation}\label{eq:quartic}
  Q_{1}(\zeta,a)\;=\;\zeta^{4}+4a^{4}\zeta^{3}+(2a^{8}+4a^{4})\zeta^{2}+4a^{4}\zeta+1 .
\end{equation}
\end{lemma}

\begin{proof}
Kernel-proved in the form \texttt{torus\_2x2\_partitionFun\_eq}; the displayed coefficients follow by reading off the sum over the sixteen spin sets, or by specializing \eqref{eq:factorization} to $n=1$ and expanding. Reciprocity is visible: the coefficient list is palindromic.
\end{proof}

\begin{lemma}[Reduction to one variable]\label{lem:reduce}
For $\zeta\ne0$, put $w=\zeta+\zeta^{-1}$. Then
\begin{equation}\label{eq:reduced}
  \frac{Q_{1}(\zeta,a)}{\zeta^{2}}\;=\;w^{2}+4a^{4}w+\bigl(4a^{4}+2a^{8}-2\bigr).
\end{equation}
\end{lemma}

\begin{proof}
Divide \eqref{eq:quartic} by $\zeta^{2}$ and group: the $\zeta^{2}+\zeta^{-2}$ term is $w^{2}-2$, and the $\zeta+\zeta^{-1}$ term is $4a^{4}w$. Kernel-proved as \texttt{recip\_quartic\_to\_quadratic}, under the hypothesis $\zeta\ne0$ (which holds for every zero, since the constant coefficient of \eqref{eq:quartic} is $1$).
\end{proof}

On the unit circle, $w=\zeta+\zeta^{-1}=2\Rea\zeta$, so the arc condition is exactly the statement that \eqref{eq:reduced} has no zero with $w>0$.

\begin{lemma}[Positivity]\label{lem:pos}
Let $b=a^{4}$ and suppose $b\ge\sqrt{2}-1$. Then for every $w>0$,
\[
  w^{2}+4bw+(4b+2b^{2}-2)\;>\;0 .
\]
\end{lemma}

\begin{proof}
Kernel-proved: \texttt{quadratic\_pos\_of\_w\_pos}. At $b=\sqrt{2}-1$ the constant term is $4(\sqrt{2}-1)+2(3-2\sqrt{2})-2=0$, so the quadratic is $w(w+4b)>0$ for $w>0$; for larger $b$ all three coefficients are nonnegative and at least one is positive.
\end{proof}

\begin{theorem}[$2\times2$ arc condition]\label{thm:2x2}
Let $a\in[0,1]$ with $a^{4}\ge\sqrt{2}-1$. Then every zero of $Q_{1}(\cdot,a)$ satisfies $\Rea\zeta\le0$.
\end{theorem}

\begin{proof}
Kernel-proved: \texttt{torus\_2x2\_arc\_condition}. Suppose $\Rea\zeta>0$. By the circle theorem $\lvert\zeta\rvert=1$, so $w=2\Rea\zeta>0$ and $\zeta\ne0$; by Lemma~\ref{lem:reduce} and $Q_{1}(\zeta,a)=0$ we get \eqref{eq:reduced}$=0$, contradicting Lemma~\ref{lem:pos}. The circle theorem enters through \texttt{isingPolynomial\_roots\_norm}, the project's Lean proof that the roots of the Ising polynomial of a ferromagnetic edge list of weight at most one lie on the unit circle.
\end{proof}

\begin{corollary}[The threshold for $n=1$]\label{cor:ac1}
The threshold of Definition~\ref{def:acn} at $n=1$ is
\[
  a_{c}(1)\;=\;\bigl(\sqrt{2}-1\bigr)^{1/4}\;=\;0.802243262923150\ldots .
\]
\end{corollary}

\begin{proof}
The hypothesis of Theorem~\ref{thm:2x2} is $a^{4}\ge\sqrt{2}-1$, i.e. $a\ge(\sqrt{2}-1)^{1/4}$ since $a\ge0$. For sharpness, substitute $\zeta=i$ in \eqref{eq:quartic}: since $i^{4}=1$, $i^{3}=-i$, $i^{2}=-1$, one gets
\[
  Q_{1}(i,a)=2-2a^{8}-4a^{4}=-2\bigl(a^{8}+2a^{4}-1\bigr),
\]
which vanishes exactly at $a^{4}=\sqrt{2}-1$. Hence at the endpoint the zero set meets the imaginary axis, and the endpoint cannot be lowered. (The sufficiency is kernel-proved; the displayed substitution is an elementary computation, and the inference that the endpoint is optimal uses the count argument of \S\ref{sec:uniform}, so this corollary is analysis-level in its sharpness half.)
\end{proof}

\begin{remark}\label{rem:2x2remark}
It is worth noting how small the threshold is in the coupling that actually appears in the Boltzmann weight. For the $2\times2$ torus every bond has weight $a^{2}$, so the threshold in the bond variable is $b=a^{2}\ge(\sqrt{2}-1)^{1/2}=0.643594\ldots$. The value $a^{4}=\sqrt{2}-1$ is the simpler way to write it because the quartic \eqref{eq:quartic} depends on $a$ only through $a^{4}$.
\end{remark}

\section{Evaluating the factorization at \texorpdfstring{$\zeta=i$}{zeta=i}}\label{sec:i}

\subsection{The cubic becomes real}\label{sec:cubicreal}

Throughout this section write
\begin{equation}\label{eq:alphabet}
  \alpha=a^{2}(1+a^{2}),\qquad \beta=1-a^{4},\qquad \gamma=a^{2}(1-a^{2})^{3},\qquad d=a^{2}(1-a^{2}),
\end{equation}
all real and, for $0<a<1$, all positive except that $\beta>0$ always.

\begin{lemma}[Cubic at $i$]\label{lem:cubici}
At $\zeta=i$ the cubic \eqref{eq:cubic} equals
\[
  X^{3}-i\alpha X^{2}-\beta X+i\gamma ,
\]
and its roots are $i y_{1},iy_{2},iy_{3}$, where $y_{1},y_{2},y_{3}$ are the roots of the real cubic
\begin{equation}\label{eq:Ra}
  R_{a}(y)\;=\;y^{3}-\alpha y^{2}+\beta y-\gamma .
\end{equation}
\end{lemma}

\begin{proof}
At $\zeta=i$ one has $\zeta^{2}=-1$ and $\zeta^{3}=-i$; substituting in \eqref{eq:tau}--\eqref{eq:delta} gives $\tau=i\alpha$, $\sigma=-\beta$, $\delta=-i\gamma$, hence the cubic $X^{3}-i\alpha X^{2}-\beta X+i\gamma$. Substituting $X=iy$ gives
\[
  (iy)^{3}-i\alpha(iy)^{2}-\beta(iy)+i\gamma
  =-iy^{3}+i\alpha y^{2}-i\beta y+i\gamma
  =-i\bigl(y^{3}-\alpha y^{2}+\beta y-\gamma\bigr),
\]
so $X$ is a root if and only if $X=iy$ with $R_{a}(y)=0$.
\end{proof}

\begin{proposition}[Reduction]\label{prop:reduction}
For every $n\ge0$ and every $a\in[0,1]$,
\begin{equation}\label{eq:Fn}
  Q_{n}(i,a)\;=\;i^{\,n+1}F_{n}(a),
  \qquad
  F_{n}(a)\;=\;p_{n+1}(a)+d(a)^{n+1},
\end{equation}
where $p_{k}=y_{1}^{k}+y_{2}^{k}+y_{3}^{k}$ is the $k$-th power sum of the roots of $R_{a}$ and $d=a^{2}(1-a^{2})$. In particular $F_{n}\in\R[a]$.
\end{proposition}

\begin{proof}
By Theorem~\ref{thm:factor}, $Q_{n}(i,a)=\Tr(\mathcal{S}(i,a)^{n+1})+(i\,d)^{n+1}$. By Lemma~\ref{lem:cubici} the eigenvalues of $\mathcal{S}(i,a)$ are $iy_{1},iy_{2},iy_{3}$, so the trace term is $\sum_{j}(iy_{j})^{n+1}=i^{n+1}p_{n+1}$. The antisymmetric term is $i^{n+1}d^{n+1}$. Adding gives \eqref{eq:Fn}. Reality is immediate from $p_{n+1},d\in\R$.
\end{proof}

\begin{remark}[Consistency with reciprocity]\label{rem:recip}
$Q_{n}$ is reciprocal of degree $2(n+1)$, because $c(S)=c(\Lambda_{n}\setminus S)$ for every $S$ and $\lvert\Lambda_{n}\rvert=2(n+1)$. Hence $Q_{n}(i)=(-1)^{n+1}\overline{Q_{n}(i)}$, which is exactly the phase statement \eqref{eq:Fn}. This is a useful check on the algebra of Lemma~\ref{lem:cubici} and on \eqref{eq:lambda4}.
\end{remark}

\begin{lemma}[Newton recurrence]\label{lem:newton}
$p_{0}=3$, $p_{1}=\alpha$, $p_{2}=\alpha^{2}-2\beta$, and for $k\ge3$,
\begin{equation}\label{eq:newton}
  p_{k}\;=\;\alpha\,p_{k-1}-\beta\,p_{k-2}+\gamma\,p_{k-3}.
\end{equation}
\end{lemma}

\begin{proof}
Newton's identities for the monic cubic \eqref{eq:Ra}.
\end{proof}

\begin{corollary}[The first two cases]\label{cor:firsttwo}
\begin{align}
  F_{1}(a)&=\;2\bigl(a^{8}+2a^{4}-1\bigr), \label{eq:F1}\\
  F_{2}(a)&=\;2a^{4}\bigl(3a^{6}+7a^{2}-6\bigr). \label{eq:F2}
\end{align}
\end{corollary}

\begin{proof}
$F_{1}=p_{2}+d^{2}=(\alpha^{2}-2\beta)+d^{2}$; expanding gives \eqref{eq:F1}. $F_{2}=p_{3}+d^{3}$ with $p_{3}=\alpha p_{2}-\beta p_{1}+\gamma\cdot 3$; expanding gives \eqref{eq:F2}.
\end{proof}

\begin{remark}\label{rem:F1again}
Corollary~\ref{cor:firsttwo} recovers Corollary~\ref{cor:ac1}: $F_{1}$ vanishes exactly at $a^{4}=\sqrt{2}-1$, so at $n=1$ the crossing happens at the $2\times2$ threshold. This is the first nontrivial check of the whole reduction, and it is exact.
\end{remark}

\subsection{Why \texorpdfstring{$F_{n}$}{Fn} decides everything}\label{sec:whyn}

Because the zeros of $Q_{n}(\cdot,a)$ lie on the unit circle, the only points of the circle with $\Rea\zeta=0$ are $\zeta=\pm i$. Since $Q_{n}$ has real coefficients, $Q_{n}(-i,a)=\overline{Q_{n}(i,a)}$; by \eqref{eq:Fn} both vanish exactly when $F_{n}(a)=0$. Thus
\[
  \text{a zero lies on the imaginary axis}
  \iff F_{n}(a)=0 .
\]
The whole threshold question for fixed $n$ is therefore the question of where the real polynomial $F_{n}$ changes sign on $(0,1)$, and the whole uniformity question is the question of whether one coupling makes all $F_{n}$ positive at once.

\section{The critical constant}\label{sec:crit}

\subsection{The discriminant}\label{sec:disc}

\begin{lemma}[Discriminant]\label{lem:disc}
For the cubic \eqref{eq:Ra},
\begin{equation}\label{eq:disc}
  \disc R_{a}\;=\;4(a^{2}-1)^{2}\,P(a^{2}),
\end{equation}
where
\begin{equation}\label{eq:octic}
  P(x)\;=\;x^{8}+2x^{7}-2x^{6}+26x^{5}-48x^{4}+30x^{3}-2x^{2}-2x-1 .
\end{equation}
\end{lemma}

\begin{proof}
For $y^{3}+by^{2}+cy+e$ the discriminant is $b^{2}c^{2}-4c^{3}-4b^{3}e-27e^{2}+18bce$; here $b=-\alpha$, $c=\beta$, $e=-\gamma$. Substituting \eqref{eq:alphabet} and collecting in $a^{2}$ gives the stated factorization. The identity has been checked symbolically against the general formula and, independently, the numerics of \S\ref{sec:numerics} evaluate both sides at the certified root of $P$.
\end{proof}

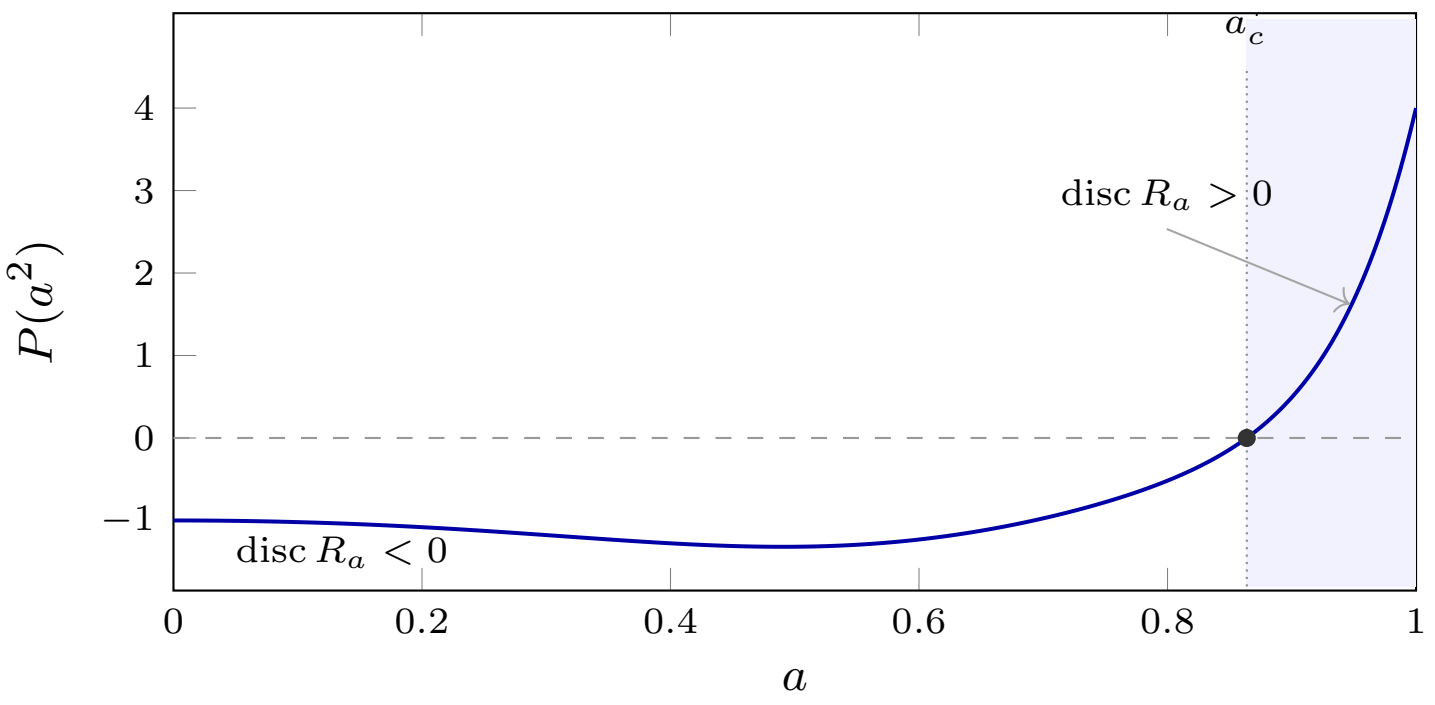
\begin{figure}[tbp]
\centering
\begin{tikzpicture}
\begin{axis}[
  width=9.8cm,height=5.4cm,
  xlabel={$a$},ylabel={$P(a^{2})$},
  xmin=0,xmax=1,ymin=-1.85,ymax=5.15,
  xtick={0,0.2,0.4,0.6,0.8,1},
  ytick={-1,0,1,2,3,4},
  tick label style={font=\footnotesize},
  label style={font=\small},
]
 \fill[blue!5] (axis cs:0.863785056687024,-1.80) rectangle (axis cs:1,5.08);
 \addplot[black!40,dashed,line width=.4pt,domain=0:1,samples=2] {0};
 \addplot[black!45,densely dotted,line width=.4pt]
   coordinates {(0.863785056687024,-1.80) (0.863785056687024,4.50)};
 \addplot[blue!65!black,line width=.75pt]
   table[x index=0,y index=1] {figures/disc.dat};
 \addplot[only marks,mark=*,mark size=1.5pt,black!80]
   coordinates {(0.863785056687024,0)};
 \node[anchor=south,font=\footnotesize] at (axis cs:0.863785,4.55)
   {$a_{c}^{*}$};
 \node[anchor=south west,font=\footnotesize] at (axis cs:0.035,-1.80)
   {$\disc R_{a}<0$};
 \node[anchor=east,font=\footnotesize] (pos) at (axis cs:0.90,2.95)
   {$\disc R_{a}>0$};
 \draw[->,black!35,line width=.4pt] (pos.south)
   -- (axis cs:0.947,1.62);
\end{axis}
\end{tikzpicture}
\caption{The sign of the discriminant, computed from the exact integer
coefficients of the octic $P$ of \eqref{eq:octic}. Since
$\disc R_{a}=4(a^{2}-1)^{2}P(a^{2})$ and $4(a^{2}-1)^{2}>0$ on $(0,1)$, the
sign of the plotted curve is the sign of $\disc R_{a}$: the cubic has three
distinct real roots where the curve is positive and one real root together
with a conjugate pair where it is negative (Remark~\ref{rem:discsign}). The
curve crosses zero once, at $a=a_{c}^{*}=0.863785056687\ldots$, the unique
root of $P$ in $(0,1)$; the shaded strip $a>a_{c}^{*}$ is the range in which
Lemma~\ref{lem:posroots} makes every $F_{n}$ positive. The data are read from
\texttt{figures/disc.dat}, generated from $P$ in exact rational arithmetic.}
\label{fig:disc}
\end{figure}

\begin{remark}\label{rem:discsign}
Since $4(a^{2}-1)^{2}\ge0$ with equality only at $a=1$, the sign of $\disc R_{a}$ on $(0,1)$ is the sign of $P(a^{2})$. A cubic with real coefficients has three distinct real roots exactly when its discriminant is positive, and one real root and a complex-conjugate pair exactly when it is negative. Figure~\ref{fig:disc} plots this sign on $(0,1)$.
\end{remark}

\subsection{The octic}\label{sec:octic}

\begin{theorem}[The critical constant]\label{thm:acstar}
The polynomial $P$ of \eqref{eq:octic} has exactly one root $x^{*}$ in $(0,1)$, and
\[
  a^{*}_{c}\;=\;\sqrt{x^{*}}\;=\;0.863785056687024029921664\ldots
\]
is the unique root in $(0,1)$ of $\disc R_{a}=0$ away from $a=1$. Moreover $P$ has exactly one further real root, at $x=-4.419423797\ldots$, and six non-real roots.
\end{theorem}

\begin{proof}[Certification]
This is a statement about the real roots of an explicit integer polynomial; it is certified by a Sturm sequence computed in exact rational arithmetic. The script \texttt{numerics\_data/torus\_ac\_threshold.py} (section~A of its output) verifies that $P$ has exactly one root in $(0,1)$, and isolates it in a rational interval whose image under $x\mapsto\sqrt{x}$ has width $1.7\times10^{-39}$. All cited digits are stable under refinement of that interval. The count of real roots is part of the same Sturm computation.
\end{proof}

\begin{remark}[What the constant means]\label{rem:meaning}
The number $a^{*}_{c}$ is not defined as an infimum or a limit; it is the unique root in $(0,1)$ of an explicit degree-eight polynomial. It is therefore an algebraic number of degree at most eight, and every digit quoted above is determined by an exact rational computation. The factor $a^{2}=x$ is the natural variable, because both the cubic \eqref{eq:Ra} and the polynomial $F_{n}$ depend on $a$ only through $a^{2}$.
\end{remark}

\begin{remark}[The superseded value]\label{rem:superseded}
An earlier draft of the formal artifact recorded the threshold for the two-row torus as $a\approx0.824$, obtained from a coarse bisection scan over circumferences $n\le13$. That value is wrong and is superseded by $a^{*}_{c}=0.863785\ldots$; the source note has been corrected. No cause for the discrepancy is asserted here: the old number was produced by a coarse scan, it is not reproducible at the certified precision, and it should not be used or cited.
\end{remark}

\section{The finite-circumference thresholds}\label{sec:finite}

\subsection{Definition and the table}\label{sec:table}

\begin{definition}\label{def:acn2}
For $n\ge1$, let $a_{c}(n)$ be the largest root in $(0,1)$ of $F_{n}$, the \emph{root threshold}. By the discussion of \S\ref{sec:whyn} and the circle theorem, $a>a_{c}(n)$ implies that no zero of $Q_{n}(\cdot,a)$ lies on the imaginary axis, while at $a=a_{c}(n)$ one does. That this root threshold is the infimum of Definition~\ref{def:acn}, that is, that the arc condition holds for every $a>a_{c}(n)$ and not merely that no zero sits at $\zeta=\pm i$, is Corollary~\ref{cor:exact}, which is analysis-level.
\end{definition}

\begin{proposition}[Basic facts]\label{prop:basic}
For every $n\ge1$: $F_{n}$ has exactly $\lceil n/2\rceil$ roots in $(0,1)$; $a_{c}(n)$ is the largest of them; $a_{c}(n)<a^{*}_{c}$; and $F_{n}(a)>0$ for every $a\in(a^{*}_{c},1)$.
\end{proposition}

\begin{proof}
The root counts and the inequality $a_{c}(n)<a^{*}_{c}$ are certified numerically (section~B of the script output). The last assertion is a consequence of Lemma~\ref{lem:posroots} below and does not depend on the numerics: for $a>a^{*}_{c}$ the cubic \eqref{eq:Ra} has three positive roots, so every term of $F_{n}$ is positive.
\end{proof}

\begin{table}[t]
\centering
\caption{The finite-circumference thresholds $a_{c}(n)$, the number of roots of $F_{n}$ in $(0,1)$, the gap to $a^{*}_{c}$, and the scaled gap $n^{2}(a^{*}_{c}-a_{c}(n))$. Values from \texttt{numerics\_data/torus\_ac\_threshold.txt}, computed at $50$ decimal digits from the exact integer coefficients of $F_{n}$.}
\label{tab:acn}
\begin{tabular}{@{}rllll@{}}
\toprule
$n$ & $a_{c}(n)$ & $\#$ roots in $(0,1)$ & $a^{*}_{c}-a_{c}(n)$ & $n^{2}\bigl(a^{*}_{c}-a_{c}(n)\bigr)$\\
\midrule
1  & 0.802243262923150152854283 & 1  & $6.1542\times10^{-2}$ & 0.06154\\
2  & 0.840357849665090902186457 & 1  & $2.3427\times10^{-2}$ & 0.09371\\
3  & 0.851222942028998351487654 & 2  & $1.2562\times10^{-2}$ & 0.11306\\
4  & 0.855900159702158341267897 & 2  & $7.8849\times10^{-3}$ & 0.12616\\
5  & 0.858361370174410782468338 & 3  & $5.4237\times10^{-3}$ & 0.13559\\
6  & 0.859821701830115537236894 & 3  & $3.9634\times10^{-3}$ & 0.14268\\
7  & 0.860760831042938025387896 & 4  & $3.0242\times10^{-3}$ & 0.14819\\
8  & 0.861400970702158984304786 & 4  & $2.3841\times10^{-3}$ & 0.15258\\
9  & 0.861857060195477934113955 & 5  & $1.9280\times10^{-3}$ & 0.15617\\
10 & 0.862193565717041416789464 & 5  & $1.5915\times10^{-3}$ & 0.15915\\
12 & 0.8626474157497453140919   & 6  & $1.1376\times10^{-3}$ & 0.16382\\
15 & 0.863035067976513456841044 & 8  & $7.4999\times10^{-4}$ & 0.16875\\
20 & 0.863350176590748232209618 & 10 & $4.3488\times10^{-4}$ & 0.17395\\
25 & 0.86350150816767391154964  & 13 & $2.8355\times10^{-4}$ & 0.17722\\
30 & 0.863585658848314806653331 & 15 & $1.9940\times10^{-4}$ & 0.17946\\
35 & 0.863637227911481085356096 & 18 & $1.4783\times10^{-4}$ & 0.18109\\
40 & 0.86367109895313283267377  & 20 & $1.1396\times10^{-4}$ & 0.18233\\
\bottomrule
\end{tabular}
\end{table}

\begin{remark}[Crowding]\label{rem:crowding}
The roots of $F_{n}$ crowd toward $a^{*}_{c}$ as $n$ grows, so the certified values in Table~\ref{tab:acn} are computed from the exact coefficients by bisection, not on a grid.
\end{remark}

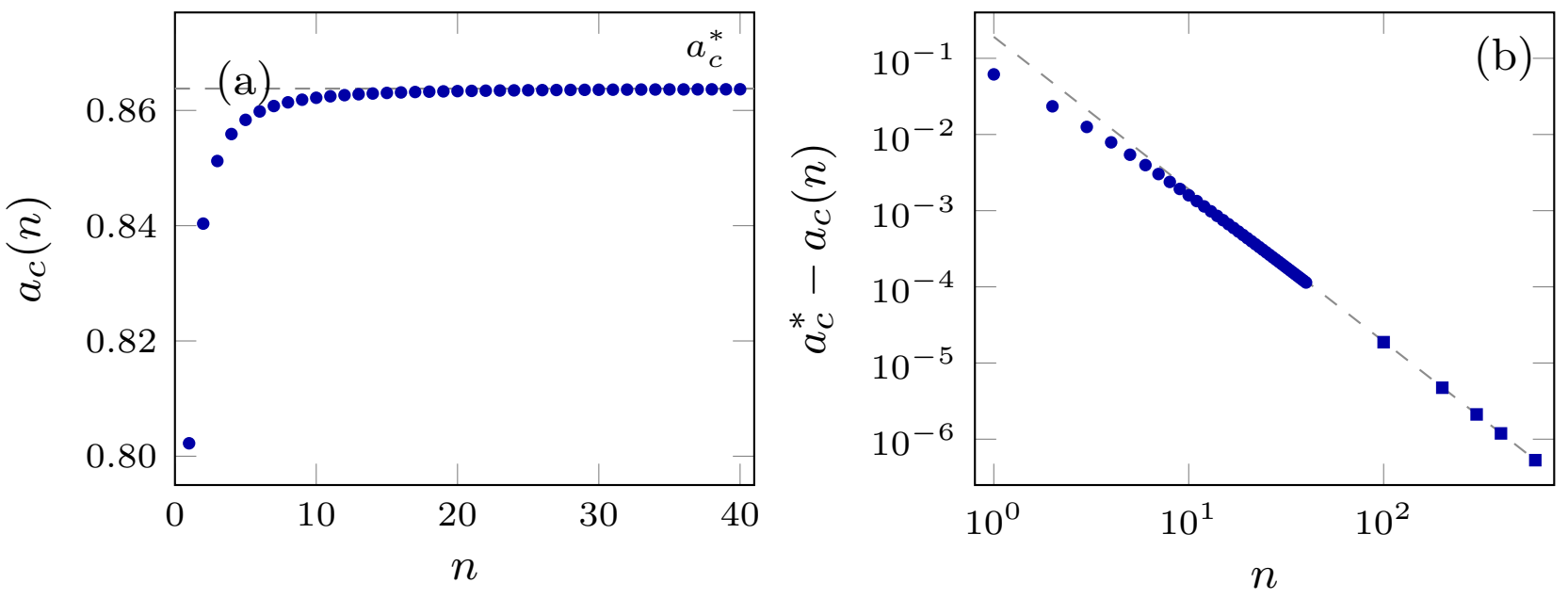
\begin{figure}[tbp]
\centering
\begin{tikzpicture}
\begin{groupplot}[
  group style={group size=2 by 1,horizontal sep=1.55cm},
  width=5.65cm,height=4.9cm,
  tick label style={font=\footnotesize},
  label style={font=\small},
]
\nextgroupplot[
  xlabel={$n$},ylabel={$a_{c}(n)$},
  xmin=0,xmax=41,ymin=0.795,ymax=0.8770,
  xtick={0,10,20,30,40},
  ytick={0.80,0.82,0.84,0.86},
  yticklabel style={/pgf/number format/fixed,
    /pgf/number format/fixed zerofill,
    /pgf/number format/precision=2},
]
 \addplot[black!45,dashed,line width=.4pt,domain=0:41,samples=2]
   {0.863785056687024};
 \addplot[only marks,mark=*,mark size=1.05pt,blue!65!black]
   table[x index=0,y index=1] {figures/acn.dat};
 \node[anchor=south east,font=\footnotesize] at (axis cs:40.4,0.86450)
   {$a_{c}^{*}$};
 \node[anchor=north west,font=\small] at (axis cs:1.6,0.8718) {(a)};
\nextgroupplot[
  xlabel={$n$},ylabel={$a_{c}^{*}-a_{c}(n)$},
  xmode=log,ymode=log,
  xmin=0.8,xmax=750,ymin=2.5e-7,ymax=4e-1,
  ytick={1e-6,1e-5,1e-4,1e-3,1e-2,1e-1},
  minor tick style={draw=none},
]
 \addplot[black!45,dashed,line width=.4pt,domain=1:600,samples=2]
   {0.191485070791/(x^2)};
 \addplot[only marks,mark=*,mark size=1.05pt,blue!65!black]
   table[x index=0,y index=2] {figures/acn.dat};
 \addplot[only marks,mark=square*,mark size=1.05pt,blue!65!black]
   table[x index=0,y index=2] {figures/acn_far.dat};
 \node[anchor=north east,font=\small] at (axis cs:740,3.4e-1) {(b)};
\end{groupplot}
\end{tikzpicture}
\caption{The finite-circumference thresholds. (a)~The certified values
$a_{c}(n)$ for $n=1,\dots,40$ (dots; Table~\ref{tab:acn} quotes a subset),
against the uniform threshold $a_{c}^{*}$ (dashed). (b)~The gap
$a_{c}^{*}-a_{c}(n)$ on a logarithmic scale, with the reference line
$C/n^{2}$; the five squares at $n=100,200,300,400,600$ are the
high-precision values of Appendix~\ref{app:lean}. Both data sets are read
from the certificate \texttt{verification/torus\_ac\_threshold.txt} (sections
B and F); nothing in the figure is typed by hand. The gap follows $C/n^{2}$
with $C=0.1914850\ldots$, the leading term of Theorem~\ref{thm:rate}.}
\label{fig:thresholds}
\end{figure}

Figure~\ref{fig:thresholds} plots the certified values against the uniform threshold, and the gap against its leading term.

\subsection{Monotonicity and convergence}\label{sec:convergence}

Numerically, $a_{c}(n)$ is strictly increasing. The rate is not an empirical observation: the discriminant zero at $a^{*}_{c}$ is simple, the double root splits with an argument $\theta$ odd in $\sqrt{a^{*}_{c}-a}$, and the threshold is the first crossing $(n+1)\theta=\pi/2$. The proof is in Appendix~\ref{app:rate}.

\begin{theorem}[Rate of convergence]\label{thm:rate}
There is $n_{0}$ such that for every $n\ge n_{0}$,
\[
  a_{c}(n)\;=\;a^{*}_{c}-\frac{C}{n^{2}}+\frac{2C}{n^{3}}-\frac{3C+D}{n^{4}}+O\bigl(n^{-5}\bigr),
  \qquad
  C=\frac{\pi^{2}\delta_{0}^{4}u^{2}}{\disc'(a^{*}_{c})}=0.1914850\ldots,
\]
where $u$ is the double root of $R_{a}$ at $a^{*}_{c}$, $y_{1}$ the simple root, $\delta_{0}=u-y_{1}$, and $D=0.130886\ldots$ is defined in Appendix~\ref{app:rate}.
\end{theorem}

The $n^{-3}$ coefficient is exact, not fitted: $R_{a}$ depends on $a$ through $a^{2}$, so the pair has even real part and odd imaginary part in $\sqrt{a^{*}_{c}-a}$, the expansion in $N=n+1$ carries no $N^{-3}$ term, and the $2C/n^{3}$ is the shift $n\to n+1$. Fitting the certified values of Appendix~\ref{app:lean} to $n=200,\dots,600$ gives $c_{1}=-0.3829\ldots$ against $-2C=-0.38297\ldots$, and $C$ itself to eight digits.

\begin{remark}[No radical formula beyond $n=1$]\label{rem:noclosed}
For $n=1$ the threshold is $(\sqrt{2}-1)^{1/4}$; for $n=2$ it is the largest root of $3a^{6}+7a^{2}-6=0$. For $n\ge3$ we claim no radical formula.
\end{remark}

\section{Uniformity in the circumference and optimality}\label{sec:uniform}

\subsection{From positive roots to the arc condition}\label{sec:posroots}

\begin{lemma}[Three positive roots give positivity]\label{lem:posroots}
Suppose $a\in(a^{*}_{c},1)$. Then $R_{a}$ has three distinct positive real roots, and consequently $F_{n}(a)>0$ for every $n\ge0$.
\end{lemma}

\begin{proof}
For $a>a^{*}_{c}$ we have $P(a^{2})>0$ (Theorem~\ref{thm:acstar} and the sign discussion of Remark~\ref{rem:discsign}), so $\disc R_{a}>0$ and $R_{a}$ has three distinct real roots. Their sum is $\alpha>0$, their product is $\gamma>0$, and the sum of pairwise products is $\beta>0$. If two of them were negative, say $-p,-q$ with $p,q>0$, and the third were $s>0$, then $\alpha>0$ forces $s>p+q$, while the pairwise sum would be $pq-s(p+q)<pq-(p+q)^{2}<0$, contradicting $\beta>0$. Hence all three are positive. Then every term of $F_{n}(a)=p_{n+1}(a)+d(a)^{n+1}$ is positive, because $p_{n+1}$ is a sum of $(n+1)$-st powers of positive numbers and $d=a^{2}(1-a^{2})>0$ for $a\in(0,1)$.
\end{proof}

\begin{theorem}[Uniform threshold]\label{thm:uniform}
Fix $a\in(a^{*}_{c},1)$. Then for every $n\ge1$, no zero of $Q_{n}(\cdot,a)$ has positive real part. Equivalently, $(a^{*}_{c},1]\subseteq A_{n}$ for every $n$.
\end{theorem}

\begin{proof}
This is an analysis-level argument. It uses the zero-count form of Rouch\'e's theorem
(Appendix~\ref{app:rouche}). The Lee--Yang circle theorem enters through the model.

\emph{Step 1.} By Lemma~\ref{lem:posroots}, $F_{n}(a)>0$, so $Q_{n}(i,a)\ne0$; by reality of the coefficients, $Q_{n}(-i,a)\ne0$ as well. Hence no zero of $Q_{n}(\cdot,a)$ lies on the imaginary axis.

\emph{Step 2.} Fix $a_{0}\in(a^{*}_{c},1)$. By Step~1, $Q_{n}(\cdot,a_{0})$ has no zero on the imaginary axis. Rouch\'e's theorem in the half-disk form gives a neighborhood of $a_{0}$ in which the number
\[
  N_{n}(a)\;=\;\#\{\zeta : Q_{n}(\zeta,a)=0,\ \Rea\zeta>0\}
\]
is constant. The interval $(a^{*}_{c},1)$ is connected, so $N_{n}$ is constant on it.

\emph{Step 3.} At $a=1$, $Q_{n}(\zeta,1)=(1+\zeta)^{2(n+1)}$, whose only zero is $\zeta=-1$, with $\Rea(-1)=-1<0$. Hence $N_{n}(1)=0$, so $N_{n}=0$ on $(a^{*}_{c},1)$.
\end{proof}

\begin{corollary}\label{cor:uniform}
For every $n\ge1$, $a_{c}(n)\le a^{*}_{c}$.
\end{corollary}

\begin{proof}
Immediate from Theorem~\ref{thm:uniform} and Definition~\ref{def:acn}.
\end{proof}

\subsection{Optimality}\label{sec:optimal}

For $a<a^{*}_{c}$ the discriminant is negative, so $R_{a}$ has one real root $y_{1}$ and a complex-conjugate pair $u\pm iv$ with $v\ne0$. Write $w=u+iv$. Then
\begin{equation}\label{eq:Fnterms}
  F_{n}(a)\;=\;y_{1}^{\,n+1}+2\Rea\bigl(w^{\,n+1}\bigr)+d^{\,n+1}.
\end{equation}

\begin{lemma}[Dominance]\label{lem:dominance}
For every $a\in(0,a^{*}_{c})$, with $y_{1}$ the real root of $R_{a}$ and $u\pm iv$ the complex pair,
\[
  \lvert u+iv\rvert\;>\;\max\bigl(y_{1},\,d(a)\bigr).
\]
\end{lemma}

\begin{proof}
Put $s=a^{2}$, so that $\alpha=s(1+s)$, $\beta=1-s^{2}$, $\gamma=s(1-s)^{3}$ and $d=s(1-s)$, and set $t=(1-s)/s$. Then
\[
  \frac{\beta}{\alpha}=\frac{1-s^{2}}{s(1+s)}=\frac{1-s}{s}=t,
  \qquad
  \frac{\gamma}{d^{2}}=\frac{s(1-s)^{3}}{s^{2}(1-s)^{2}}=\frac{1-s}{s}=t :
\]
the two bounds in the statement are the same number. Since $R_{a}(0)=-\gamma<0$ and
\[
  R_{a}(t)=\frac{\beta^{3}}{\alpha^{3}}-\gamma=\frac{(1-s)^{4}(1+s)(1+s^{2})}{s^{3}}>0,
\]
and since $\disc R_{a}<0$ on $(0,a^{*}_{c})$ (Remark~\ref{rem:discsign} and Theorem~\ref{thm:acstar}), $R_{a}$ has exactly one real root and it lies in $(0,t)$; hence $0<y_{1}<t$.

By Vieta, $2u=\alpha-y_{1}$ and $\lvert w\rvert^{2}=\beta-\alpha y_{1}+y_{1}^{2}$ with $w=u+iv$, while $y_{1}\lvert w\rvert^{2}=\gamma$. Therefore
\[
  \lvert w\rvert^{2}-y_{1}^{2}=\beta-\alpha y_{1}=\alpha\Bigl(\frac{\beta}{\alpha}-y_{1}\Bigr)=\alpha\,(t-y_{1})>0,
\]
and
\[
  \lvert w\rvert^{2}-d^{2}=\frac{\gamma}{y_{1}}-d^{2}=d^{2}\Bigl(\frac{t}{y_{1}}-1\Bigr)>0 .
\]
Hence $\lvert w\rvert>y_{1}$ and $\lvert w\rvert>d$.
\end{proof}

\begin{table}[t]
\centering
\caption{The dominance inequality of Lemma~\ref{lem:dominance}, sampled as a check. The column ``min $n$'' is the least $n$ with $F_{n}(a)<0$; it agrees with the thresholds of Table~\ref{tab:acn}.}
\label{tab:dominance}
\begin{tabular}{@{}lrrrr@{}}
\toprule
$a$ & $y_{1}$ & $\lvert u+iv\rvert$ & $d$ & min $n$\\
\midrule
0.30   & 0.06851524 & 0.99492361 & 0.08190 & 1\\
0.50   & 0.11529624 & 0.95643251 & 0.18750 & 1\\
0.60   & 0.11403452 & 0.90971016 & 0.23040 & 1\\
0.70   & 0.09275170 & 0.83712894 & 0.24990 & 1\\
0.80   & 0.05582042 & 0.73138690 & 0.23040 & 1\\
0.85   & 0.03548512 & 0.65961466 & 0.20049 & 3\\
0.8637 & 0.03018296 & 0.63648738 & 0.18949 & 47\\
\bottomrule
\end{tabular}
\end{table}

\begin{proposition}[Optimality]\label{prop:optimal}
For every $a\in(0,a^{*}_{c})$ there is $n$ with $a_{c}(n)>a$. Consequently
\[
  \sup_{n\ge1}a_{c}(n)\;=\;a^{*}_{c}.
\]
\end{proposition}

\begin{proof}
This is analysis-level. Since $\lvert w\rvert>\max(y_{1},d)$, the middle term of \eqref{eq:Fnterms} dominates: writing $w=\lvert w\rvert e^{i\varphi}$,
\[
  F_{n}(a)=\lvert w\rvert^{n+1}\Bigl(2\cos\bigl((n+1)\varphi\bigr)+\Bigl(\tfrac{y_{1}}{\lvert w\rvert}\Bigr)^{n+1}+\Bigl(\tfrac{d}{\lvert w\rvert}\Bigr)^{n+1}\Bigr),
\]
and the last two summands tend to $0$, so for all large $n$ the sign of $F_{n}(a)$ is the sign of $\cos((n+1)\varphi)$ whenever that cosine is bounded away from $0$. As $n$ ranges over the integers, $(n+1)\varphi$ is equidistributed modulo $2\pi$ when $\varphi/2\pi$ is irrational, and periodic with some period $q\ge3$ when it is rational (the cases $q=1,2$ are excluded by $v\ne0$); in either case there are infinitely many $n$ with $\cos((n+1)\varphi)\le-\tfrac12$, and then $F_{n}(a)<0$ for all sufficiently large such $n$. Fix one.

Now $F_{n}(a^{*}_{c})>0$ by Lemma~\ref{lem:posroots} (applied at any $a'>a^{*}_{c}$ and extended to $a^{*}_{c}$ by continuity, or directly: at $a^{*}_{c}$ the roots of $R_{a}$ are real and nonnegative, and $\gamma>0$ at $a^{*}_{c}$ forces them positive). Since $F_{n}(a)<0$ and $F_{n}$ is continuous, $F_{n}$ has a root $a'\in(a,a^{*}_{c})$. At $a'$ we have $Q_{n}(i,a')=0$. By the counting argument of Theorem~\ref{thm:uniform}, $N_{n}$ is constant on each interval between consecutive roots of $F_{n}$ in $(a,a^{*}_{c})$, and it equals $0$ on $(a',a^{*}_{c})$. At the crossing $a'$ a zero passes through $i$, so on the interval immediately below $a'$ the count $N_{n}$ is at least $1$. Therefore $Q_{n}(\cdot,a'')$ fails the arc condition for some $a''\in(a,a')$, i.e.\ $a_{c}(n)>a$. Since $a<a^{*}_{c}$ was arbitrary, $\sup_{n}a_{c}(n)\ge a^{*}_{c}$; combined with Corollary~\ref{cor:uniform}, the supremum equals $a^{*}_{c}$.
\end{proof}

\begin{corollary}[The uniform threshold is exactly $a^{*}_{c}$]\label{cor:exact}
The arc condition holds for $Q_{n}(\cdot,a)$ for every $n\ge1$ if and only if $a\ge a^{*}_{c}$.
\end{corollary}

\begin{proof}
Sufficiency is Theorem~\ref{thm:uniform}; necessity is Proposition~\ref{prop:optimal}, which shows that every $a<a^{*}_{c}$ fails for some $n$.
\end{proof}

\begin{remark}[The status of Corollary~\ref{cor:exact}]\label{rem:statusexact}
Corollary~\ref{cor:exact} is the cleanest statement of the paper's main mathematical conclusion. Both directions are analysis-level: sufficiency rests on the zero-count form of Rouch\'e's theorem (Appendix~\ref{app:rouche}); necessity on the dominance lemma, proved above. The remaining formalization gap is the algebraic $\zeta=i$ reduction and the analytic implicit function theorem, not the counting step. A reader who wants only formally checked statements should read \S\ref{sec:factor}--\S\ref{sec:two} and Appendix~\ref{app:lean}; a reader who wants the sharp constant should read the whole paper.
\end{remark}

\section{Discussion and outlook}\label{sec:discussion}

\subsection{Two directions of exact solvability}\label{sec:twodirections}

The two-dimensional Ising model is exactly solvable, and it is worth being precise about which exact solvability is used here. Onsager's solution \cite{Onsager1944} gives the free energy in zero field; the transfer-matrix argument was recast by Kaufman \cite{Kaufman1949} as a spinor diagonalization of the row-to-row matrix of a strip of height $h$, a matrix of size $2^h$, on the finite periodic lattice, and that analysis is given a representation-theoretic treatment by Hystad \cite{Hystad2011}. The thermodynamics of bounded strips, including the finite-circumference corrections, is the subject of Ferdinand and Fisher \cite{FerdinandFisher1969}.

All of this is zero-field exact solvability, and the restriction is not incidental. The spinor transformation requires a conserved spin-flip parity, which splits the transfer matrix into sectors and turns it into a free-fermion operator; a magnetic field breaks that parity, and with it the free-fermion description. Here the field is the site fugacity $\zeta$: $Q_n(\zeta,a)$ is an Ising partition function in a field, with $\zeta$ carried as a variable rather than fixed at a number. The classical analysis is therefore unavailable not because the strip is finite but because the field is present.

What replaces it is a different symmetry. At height two the column has four configurations and the transfer matrix has size $4$; the row interchange commutes with it for every $\zeta$; and the interchange splits the space into a one-dimensional antisymmetric line and a three-dimensional symmetric block. The antisymmetric eigenvalue is the monomial $\zeta a^2(1-a^2)$, and the symmetric block has an explicit cubic. The circumference is arbitrary and finite and no limit is taken, so the reward is a closed form valid at every circumference.

\subsection{Why uniformity in the circumference matters}\label{sec:whyuniform}

The arc condition is a hypothesis in a local limit theorem, and a hypothesis that has to be retuned as the system grows is useless for a limit theorem. Suppose one wants a local limit theorem for the occupation count on the sequence of tori $\Lambda_n$ at a fixed coupling $a$. The Lee--Yang route delivers a bound of order $\lvert\Lambda_n\rvert^{-1}$ provided the arc condition holds at that $a$ for each $n$. If the threshold grew with $n$, say $a_c(n)\to1$, then no fixed $a$ would work and the route would give nothing in the limit. Theorem~\ref{thm:uniform} and Corollary~\ref{cor:exact} say the opposite: the thresholds are bounded above by $a_c^*<1$, so every $a\in(a_c^*,1)$ works for all $n$ at once.

The rate is the quantitative form of the same statement. It is not enough that $a_c(n)$ stay below a fixed constant; one wants to know how fast the finite thresholds reach it, because that controls the coupling margin in a finite-circumference application. Theorem~\ref{thm:rate} gives the expansion with an explicit algebraic leading constant.

\subsection{Formalization}\label{sec:roleformal}

The kernel-proved part of the paper is not a re-derivation of known facts: the factorization, the cubic with its exact coefficients, the Cayley--Hamilton identity and the recurrence are specific to this model, and the $2\times2$ threshold is a sharp constant. Formalizing them has two benefits. First, the identities are checked against the actual transfer matrix, not a hand-simplified surrogate. Second, the division of labor is explicit: the artifact certifies the algebra; the analysis-level layer is confined to the classical inputs and the local expansion of Appendix~\ref{app:rate}.

The counting input of Theorem~\ref{thm:uniform} is formalized; Appendix~\ref{app:rouche} records the two proofs and the comparator output. The dominance lemma is proved in \S\ref{sec:uniform}; it is analysis-level, and this paper does not claim a Lean proof of it. The algebraic $\zeta=i$ reduction, the uniformity theorem and the rate expansion are analysis-level, not Lean. These are formalization tasks, not open mathematical problems.

\subsection{Relation to the literature}\label{sec:lit}

The circle theorem is due to Lee and Yang \cite{YangLee1952,LeeYang1952}, with the contraction argument of Asano \cite{Asano1970}, the extension of Ruelle \cite{Ruelle1971}, the generalized-systems version of Newman \cite{Newman1974}, and the general multicomponent version of Lieb and Sokal \cite{LiebSokal1981}; a modern account of the analytic consequences, including limit theorems, is \cite{LPRS2016}. The transfer-matrix formalism for strips is classical and is treated in \cite{Baxter1982}. The Ising model goes back to \cite{Ising1925} and its two-dimensional solution to \cite{Onsager1944}. The motion of Lee--Yang zeros as a coupling varies is itself studied: for ferromagnetic Ising systems in a complex field the zeros move monotonically and distinct trajectories do not meet \cite{HouJiangNewman2023}; here the varying parameter is the coupling $a$, and the analysis-level argument uses continuity and the equimodular limit \cite{BerahaKahaneWeiss1975,SalasSokal2001}.

What appears to be new is the combination at this scale: the exact reflection decomposition of the two-row transfer matrix, the closed form uniform in the circumference, the explicit algebraic $a_c^*$ and the proof that it is exactly the uniform threshold, and the second-order rate at which the finite thresholds approach it.

\subsection{Outlook}\label{sec:open}

The two-row threshold problem is complete at the analysis level. Three directions lie outside this paper.

\begin{enumerate}
\item \emph{Height three and beyond.} The row interchange still splits the transfer matrix, but the symmetric block is larger and its factorization changes. A characteristic polynomial uniform in the height, and the $h\to\infty$ limit, are open.
\item \emph{Fixed boundary conditions.} The circle theorem fails for fixed-end cylinders. Exact small-width algebra and numerical crossings do not constitute an all-width theorem.
\item \emph{From the arc condition to a limit theorem.} The present paper supplies the threshold and its rate; turning the arc condition into a sharp local limit theorem at the threshold is a separate analytic step, outside the scope of this paper.
\end{enumerate}

None of these affects the two-row results stated here.

\appendix
\section{The formal artifact}\label{app:lean}

\subsection{Modules and build}\label{app:modules}

The formalization lives in two modules of the Lean~4 / Mathlib project \texttt{SteinGate} \cite{deMoura2021,mathlib2020}:
\begin{itemize}
\item \texttt{lean/SteinGate/LeeYang/LeeYang2DTorus.lean}: the trace identity on a genuine $(m+1)\times(n+1)$ periodic torus;
\item \texttt{lean/SteinGate/LeeYang/LeeYang2DTorusArcCondition.lean}: the $2\times(n+1)$ specialization: entry lemmas, reflection decomposition, factorization, cubic, Cayley--Hamilton, recurrence, and the $2\times2$ threshold.
\end{itemize}
The circle theorem at every torus is instantiated in the first module as
\texttt{torusEdges\_roots\_norm} (all $(m+1)\times(n+1)$ tori, $0\le a\le1$). From the repository root,
\begin{verbatim}
lake build SteinGate.LeeYang.LeeYang2DTorus
lake build SteinGate.LeeYang.LeeYang2DTorusArcCondition
\end{verbatim}
Both modules build with exit code $0$. Neither contains a \texttt{sorry}. The axiom footprint of every headline declaration is \texttt{[propext, Classical.choice, Quot.sound]}, i.e.\ Mathlib's own; the check is
\begin{verbatim}
#print axioms SteinGate.isingPartitionFun_torusEdges_2xn_eq
#print axioms SteinGate.charpoly_torusT_sym
#print axioms SteinGate.torusT_symTrace_recurrence
#print axioms SteinGate.torus_2x2_arc_condition
\end{verbatim}
executed in a file that imports the two modules.

\subsection{Declaration map}\label{app:decl}

Table~\ref{tab:decl} lists every declaration the paper leans on, the line at which it occurs in
the current revision, and the statement it certifies.

\begin{table}[h]
\centering
\caption{The declarations used in this paper, with the statement each certifies. Line numbers refer to the current revision of the two modules.}
\label{tab:decl}
\footnotesize\setlength{\tabcolsep}{4pt}
\begin{tabular}{@{}p{0.40\textwidth}lp{0.46\textwidth}@{}}
\toprule
Declaration & Line & Statement \\
\midrule
\multicolumn{3}{@{}l}{\emph{\texttt{LeeYang2DTorus.lean}}}\\
\texttt{torusEdges\_\allowbreak{}weight\_\allowbreak{}eq} & 114 & Every edge of the torus carries weight $a$.\\
\texttt{torusEdges\_\allowbreak{}roots\_\allowbreak{}norm} & 126 & Circle theorem at every $(m+1)\times(n+1)$ torus, $0\le a\le1$.\\
\texttt{isingPartitionFun\_\allowbreak{}torusEdges\_\allowbreak{}eq\_\allowbreak{}trace} & 293 & $\displaystyle Z_{m,n}(\zeta,a)=\Tr(\mathcal{T}_{m}^{n+1})$ for the $(m+1)\times(n+1)$ torus.\\
\addlinespace
\multicolumn{3}{@{}l}{\emph{\texttt{LeeYang2DTorusArcCondition.lean}}}\\
\texttt{T\_\allowbreak{}c00\_\allowbreak{}c00} \dots\ \texttt{T\_\allowbreak{}c11\_\allowbreak{}c11} & 77--122 & The sixteen entries of \eqref{eq:T}.\\
\texttt{torus\_\allowbreak{}2x2\_\allowbreak{}partitionFun\_\allowbreak{}eq} & 170 & The quartic \eqref{eq:quartic}.\\
\texttt{recip\_\allowbreak{}quartic\_\allowbreak{}to\_\allowbreak{}quadratic} & 207 & The reduction \eqref{eq:reduced} for $\zeta\ne0$.\\
\texttt{quadratic\_\allowbreak{}pos\_\allowbreak{}of\_\allowbreak{}w\_\allowbreak{}pos} & 237 & Lemma~\ref{lem:pos}.\\
\texttt{torus\_\allowbreak{}2x2\_\allowbreak{}arc\_\allowbreak{}condition} & 273 & Theorem~\ref{thm:2x2}.\\
\texttt{torusT\_\allowbreak{}antisym\_\allowbreak{}eigenvalue} & 380 & Lemma~\ref{lem:antisym}.\\
\texttt{torusT\_\allowbreak{}sym\_\allowbreak{}charPoly} & 434 & The cubic \eqref{eq:cubic}, with the three coefficients named.\\
\texttt{torusT\_\allowbreak{}sym} & 458 & The block \eqref{eq:S}.\\
\texttt{torusT\_\allowbreak{}mulVec\_\allowbreak{}symVec} & 525 & Lemma~\ref{lem:invariance}.\\
\texttt{trace\_\allowbreak{}pow\_\allowbreak{}torusT} & 635 & Theorem~\ref{thm:split}.\\
\texttt{isingPartitionFun\_\allowbreak{}torusEdges\_\allowbreak{}2xn\_\allowbreak{}eq} & 648 & Theorem~\ref{thm:factor}.\\
\texttt{charpoly\_\allowbreak{}torusT\_\allowbreak{}sym} & 658 & Theorem~\ref{thm:cubic}: the block's characteristic polynomial is the cubic.\\
\texttt{torusT\_\allowbreak{}sym\_\allowbreak{}charPoly\_\allowbreak{}eq} & 689 & The cubic in expanded form, with $\tau,\sigma,\delta$ of \eqref{eq:tau}--\eqref{eq:delta}.\\
\texttt{torusT\_\allowbreak{}sym\_\allowbreak{}cayley\_\allowbreak{}hamilton} & 697 & Theorem~\ref{thm:CH}.\\
\texttt{torusT\_\allowbreak{}symTrace\_\allowbreak{}recurrence} & 713 & Theorem~\ref{thm:rec}.\\
\bottomrule
\end{tabular}
\end{table}

\subsection{Numerics}\label{sec:numerics}

The computations behind \S\ref{sec:crit}--\S\ref{sec:uniform} are produced by
\begin{verbatim}
python3 verification/torus_ac_threshold.py
\end{verbatim}
which writes \texttt{verification/torus\_ac\_threshold.txt}. Dependencies are \texttt{mpmath} and \texttt{sympy}; the runtime is about half a minute. Section~A of the output certifies $x^{*}$ and $a^{*}_{c}$ by Sturm sequences in exact rational arithmetic (bracket width $1.7\times10^{-39}$). Section~B gives $a_{c}(n)$ for $n=1,\dots,40$ at $24$ significant digits, from the exact integer coefficients of $F_{n}$ and with no floating-point cancellation. Section~C checks \eqref{eq:Fn} against a direct evaluation of the $3\times3$ block at $50$ digits, with maximum discrepancy $1.022\times10^{-48}$. Section~D checks that $R_{a}$ has three positive roots on a $401$-point grid over $(a^{*}_{c},1)$, with no failures. Section~E re-checks the dominance inequality of Lemma~\ref{lem:dominance} on $2000$ interior points of $(0,a^{*}_{c})$, with no failures, now a consistency check since the inequality is proved in \S\ref{sec:optimal}, and reports the least $n$ with $F_{n}(a)<0$ at seven sample points. Section~F computes $a_{c}(n)$ for $n=100,200,300,400,600$ at $70$ digits from the power-sum recurrence, fits the constant of Theorem~\ref{thm:rate}, and checks the mechanism identity $\cos((n+1)\theta)=-\varepsilon_{n}$ at $n=10,20,40$. The cylinder thresholds of \S\ref{sec:open} are produced by \texttt{verification/cylinder\_boundaries.py}.

The digits quoted for $a^{*}_{c}$ are certified, not merely computed: the Sturm bracket isolating it has width $1.7\times10^{-39}$, so all $24$ digits displayed are determined, and the trailing ``$\ldots$'' indicates that the expansion continues.

\subsection{The Rouch\'e counting input}\label{app:rouche}

The zero-count form of Rouch\'e's theorem is formalized outside the two modules above. The
LeanEval problem \texttt{rouche\_zero\_count\_eq} asks for equality of the multiplicity-counted
zero counts of $f$ and $f+g$ on a disk under the boundary bound. Two proofs pass the
official comparator at Lean \texttt{v4.32.2} and mathlib \texttt{905b9581}: a port of the disk
theorem used in \S\ref{sec:uniform}, and a direct construction of the argument principle from the
divisor API. The sources are recorded at commits \texttt{5c440946} and \texttt{cc0f06d} in the
development repository; they are available from the author while the repository is private. The comparator ran the statement comparison, the axiom check, the Lean
kernel, and nanoda; both proofs end with \texttt{Your solution is okay!}. The run was on macOS
with the repository's no-op \texttt{landrun} shim; the statement, axiom, Lean kernel, and nanoda
checks ran.

The theorem supplies the counting input in Theorem~\ref{thm:uniform}. The half-disk form used
there follows by the standard contour argument. The remaining Lean work is the algebraic
$\zeta=i$ reduction and a definition of $a^{*}_{c}$; the analytic implicit function theorem for
Theorem~\ref{thm:rate} is a separate port. The paper's artifact does not import this development.

\section{Proof of Theorem~\ref{thm:rate}}\label{app:rate}

Throughout, $\delta=a^{*}_{c}-a$, $s=\sqrt{\delta}$, $N=n+1$, and $\alpha,\beta,\gamma,d$ are as in \eqref{eq:alphabet} and \S\ref{sec:i}. For $a<a^{*}_{c}$ near $a^{*}_{c}$ write the roots of $R_{a}$ as $y_{1}(a)$ and $u(a)\pm iv(a)$ with $v>0$, and put $u_{\ast}:=u(a^{*}_{c})$, $y_{1\ast}:=y_{1}(a^{*}_{c})$, $\delta_{0}=u_{\ast}-y_{1\ast}>0$, $r=\lvert u+iv\rvert$.

\begin{lemma}[Splitting]\label{lem:split}
There are $\eta>0$ and real-analytic functions $U,W$ on $[0,\eta^{2})$ such that
\[
  u(a)=U(\delta),\qquad v(a)=s\,W(\delta),\qquad U(0)=u_{\ast},\qquad
  W(0)=c:=\frac{\sqrt{\disc'(a^{*}_{c})}}{2\delta_{0}^{2}}>0 .
\]
\end{lemma}

\begin{proof}
$y_{1\ast}$ is a simple root of $R_{a^{*}_{c}}$, so by the analytic implicit function theorem $y_{1}$ extends to a real-analytic function of $a$ near $a^{*}_{c}$. By Vieta $2u=\alpha-y_{1}$, so $u=U(\delta)$ with $U$ real-analytic and $U(0)=u_{\ast}$.

Again by Vieta, $\lvert w\rvert^{2}=\beta-\alpha y_{1}+y_{1}^{2}$; hence $v^{2}=\lvert w\rvert^{2}-u^{2}$ is real-analytic in $a$ and vanishes at $a^{*}_{c}$. The discriminant identity $\disc R_{a}=-4v^{2}\bigl[(y_{1}-u)^{2}+v^{2}\bigr]^{2}$ gives at $a=a^{*}_{c}$
\[
  \disc'(a^{*}_{c})=-4\delta_{0}^{4}\,(v^{2})'(a^{*}_{c}),
\]
so $v^{2}=(a^{*}_{c}-a)W_{0}(a)$ with $W_{0}$ real-analytic and $W_{0}(a^{*}_{c})=\disc'(a^{*}_{c})/(4\delta_{0}^{4})>0$. Since $W_{0}>0$ near $a^{*}_{c}$, $W(\delta):=\sqrt{W_{0}(a^{*}_{c}-\delta)}$ is real-analytic and $v=sW(\delta)$.
\end{proof}

\begin{lemma}[Odd angle]\label{lem:odd}
Read as a function of $s$, $\theta(a)=\arctan\bigl(v(a)/u(a)\bigr)$ is real-analytic and odd on $[0,s_{0})$ for some $s_{0}>0$, and
\[
  \theta(s)=\kappa s+\kappa_{3}s^{3}+O(s^{5}),\qquad \kappa=\frac{c}{u_{\ast}}>0 .
\]
\end{lemma}

\begin{proof}
By Lemma~\ref{lem:split}, $u=U(s^{2})$ is even and $v=sW(s^{2})$ is odd, with $U(0)=u_{\ast}>0$; so $v/u$ is odd and real-analytic near $0$, and $\arctan$ preserves oddness. Since $\theta$ is odd and real-analytic, only odd powers occur; differentiating gives $\theta'(0)=W(0)/U(0)=c/u_{\ast}$, and $\kappa_{3}=\theta'''(0)/6$.
\end{proof}

\begin{lemma}[First crossing]\label{lem:first}
There is $\eta>0$ such that $r(a)>\max(y_{1}(a),d(a))$ on $[a^{*}_{c}-\eta,a^{*}_{c}]$ and $\theta$ is strictly decreasing in $a$ there. For all $n$ large enough, the largest zero of $F_{n}$ in $(a^{*}_{c}-\eta,a^{*}_{c})$ is the unique solution of
\[
  N\theta(a)=\frac{\pi}{2}+\arcsin\varepsilon_{n}(a),
  \qquad
  \varepsilon_{n}(a)=\frac{y_{1}(a)^{N}+d(a)^{N}}{2r(a)^{N}} .
\]
\end{lemma}

\begin{proof}
The inequality $r>\max(y_{1},d)$ is Lemma~\ref{lem:dominance} for $a<a^{*}_{c}$ and holds at $a^{*}_{c}$ by continuity, since $y_{1\ast}=0.030\ldots$, $d(a^{*}_{c})=0.189\ldots$ and $r(a^{*}_{c})=u_{\ast}=0.636\ldots$. By Lemma~\ref{lem:odd}, $\theta=\kappa s+O(s^{3})$ is strictly increasing in $s$, hence strictly decreasing in $a$, near $a^{*}_{c}$.

By Vieta and De Moivre,
\[
  F_{n}=y_{1}^{N}+2r^{N}\cos(N\theta)+d^{N}=2r^{N}\bigl(\cos(N\theta)+\varepsilon_{n}\bigr),
\]
so $F_{n}=0$ exactly when $\cos(N\theta)=-\varepsilon_{n}$. On the branch $N\theta\in(\pi/2,\pi)$ the unique solution is $N\theta=\pi/2+\arcsin\varepsilon_{n}$. For larger $a$, $\theta$ is smaller, so $N\theta<\pi/2+\arcsin\varepsilon_{n}$, hence $\cos(N\theta)>-\varepsilon_{n}$ and $F_{n}>0$. Thus this crossing is the largest zero in the interval.
\end{proof}

\begin{lemma}[Inversion]\label{lem:inv}
Let $s_{n}$ be defined by $N\theta(a^{*}_{c}-s_{n}^{2})=\pi/2+\arcsin\varepsilon_{n}(a^{*}_{c}-s_{n}^{2})$. Then
\[
  s_{n}=\frac{a}{N}+\frac{b}{N^{3}}+O(N^{-5}),\qquad a=\frac{\pi}{2\kappa},\quad b=-\frac{\kappa_{3}a^{3}}{\kappa}.
\]
\end{lemma}

\begin{proof}
By Lemma~\ref{lem:odd}, $\theta(s)=\kappa s+\kappa_{3}s^{3}+O(s^{5})$ with a remainder uniform on $[0,s_{0})$. Since $y_{1}/r$ and $d/r$ are bounded away from $1$ uniformly in $a$, there is $c'>0$ with $\varepsilon_{n}\le e^{-c'N}$, and $\lvert\varepsilon_{n}'\rvert\le C'N e^{-c'N}$. Put $H(s)=N\theta(s)-\pi/2-\arcsin\varepsilon_{n}(a^{*}_{c}-s^{2})$. At $s=a/N+b/N^{3}$ one has $H=\kappa a-\pi/2+O(N^{-4})$, and $H'(s)\ge N\kappa/2$ for small $s$; hence the zero of $H$ is $a/N+b/N^{3}+O(N^{-5})$.
\end{proof}

\begin{proof}[Proof of Theorem~\ref{thm:rate}]
With $s_{n}$ as in Lemma~\ref{lem:inv}, Lemma~\ref{lem:first} gives $a_{c}(n)=a^{*}_{c}-s_{n}^{2}$ for all large $n$. Hence
\[
  a_{c}(n)=a^{*}_{c}-\frac{a^{2}}{N^{2}}-\frac{2ab}{N^{4}}+O(N^{-6})
  =a^{*}_{c}-\frac{C}{N^{2}}-\frac{D}{N^{4}}+O(N^{-6}),
\]
with
\[
  C=a^{2}=\frac{\pi^{2}}{4\kappa^{2}}=\frac{\pi^{2}u_{\ast}^{2}}{4c^{2}}=\frac{\pi^{2}\delta_{0}^{4}u_{\ast}^{2}}{\disc'(a^{*}_{c})},
  \qquad
  D=2ab=-\frac{2\kappa_{3}a^{4}}{\kappa},
\]
using $c^{2}=\disc'(a^{*}_{c})/(4\delta_{0}^{4})$; numerically $\kappa_{3}=-6.406867\ldots$ and $D=0.130886\ldots$. Finally $N^{-2}=n^{-2}-2n^{-3}+3n^{-4}+O(n^{-5})$ and $N^{-4}=n^{-4}+O(n^{-5})$, which gives the stated form with the coefficient $3C+D$.
\end{proof}

\section*{Acknowledgments}

The author thanks the maintainers of Mathlib for the library on which the formalization rests.

\section*{Statements and Declarations}

\textbf{Funding.} This research received no external funding. It was supported by the author's personal funds.

\textbf{Competing interests.} The author declares no competing interests.

\textbf{Author contributions.} The author is the sole author and is responsible for the conception of the work, the mathematical arguments, the computations and the formalization.

\textbf{Use of AI tools.} Large language models (Claude, Anthropic) were used during the preparation of this work beyond copy editing: to restructure and draft parts of the text, to write and run computer-algebra and numerical scripts, and to help develop and check parts of the Lean~4 proofs. The Lean results were checked by the Lean kernel, and the exact computations were reproduced with independent scripts. The author takes full responsibility for the content of the paper.

\textbf{Ethics approval.} Not applicable; the work involves no human participants or animals.

\textbf{Data availability.} The frozen numerical output underlying the tables is \texttt{numerics\_data/torus\_ac\_threshold.txt}, regenerated by \texttt{torus\_ac\_threshold.py}. The development repository is \url{https://github.com/mrphilo420/Project-SteinGate}; it is private until publication, and the data and scripts are available from the author on reasonable request.

\textbf{Code availability.} The Lean~4 modules and the Python scripts are part of the same repository, together with instructions for reproducing the Lean results. They are available from the author on reasonable request while the repository is private.

\bibliography{references}

\end{document}